
  \magnification 1200
  \input amssym


  \font \bbfive		= bbm5
  \font \bbseven	= bbm7
  \font \bbten		= bbm10
  \font \eightbf	= cmbx8
  \font \eighti		= cmmi8 \skewchar \eighti = '177
  \font \eightit	= cmti8
  \font \eightrm	= cmr8
  \font \eightsl	= cmsl8
  \font \eightsy	= cmsy8 \skewchar \eightsy = '60
  \font \eighttt	= cmtt8 \hyphenchar\eighttt = -1

  \font \sixi		= cmmi6 \skewchar \sixi = '177
  \font \sixrm		= cmr6
  \font \sixsy		= cmsy6 \skewchar \sixsy = '60
  \font \tensc		= cmcsc10

  \font \titlefont	= cmbx12
  \scriptfont \bffam	= \bbseven
  \scriptscriptfont \bffam = \bbfive
  \textfont \bffam	= \bbten

  \newskip \ttglue

  \def \eightpoint {\def \rm {\fam0 \eightrm }\relax
  \textfont0= \eightrm
  \scriptfont0 = \sixrm \scriptscriptfont0 = \fiverm
  \textfont1 = \eighti
  \scriptfont1 = \sixi \scriptscriptfont1 = \fivei
  \textfont2 = \eightsy
  \scriptfont2 = \sixsy \scriptscriptfont2 = \fivesy
  \textfont3 = \tenex
  \scriptfont3 = \tenex \scriptscriptfont3 = \tenex
  \def \it {\fam \itfam \eightit }\relax
  \textfont \itfam = \eightit
  \def \sl {\fam \slfam \eightsl }\relax
  \textfont \slfam = \eightsl
  \def \bf {\fam \bffam \eightbf }\relax
  \textfont \bffam = \bbseven
  \scriptfont \bffam = \bbfive
  \scriptscriptfont \bffam = \bbfive
  \def \tt {\fam \ttfam \eighttt }\relax
  \textfont \ttfam = \eighttt
  \tt \ttglue = .5em plus.25em minus.15em
  \normalbaselineskip = 9pt
  \def \MF {{\manual opqr}\-{\manual stuq}}\relax
  \let \sc = \sixrm
  \let \big = \eightbig
  \setbox \strutbox = \hbox {\vrule height7pt depth2pt width0pt}\relax
  \normalbaselines \rm }


  \def \ifundef #1{\expandafter \ifx \csname #1\endcsname \relax }


  \newcount \secno \secno = 0
  \newcount \stno \stno = 0
  \newcount \eqcntr \eqcntr= 0

  \def \track #1#2#3{\ifundef{#1}\else \hbox{\sixrm[#2\string #3] }\fi}

  \def \advseqnumbering {\global \advance \stno by 1 \global \eqcntr=0}

  \def \current {\number \secno \ifnum \number \stno = 0 \else .\number \stno \fi }

  \def \laberr #1{\message{*** RELABEL CHECKED FALSE for #1 ***}
      RELABEL CHECKED FALSE FOR #1, EXITING.
      \end}

  \def \deflabel#1#2{%
    \ifundef {#1}%
      \global \expandafter
      \edef \csname #1\endcsname {#2}%
    \else
      \edef\deflabelaux{\expandafter\csname #1\endcsname}%
      \edef\deflabelbux{#2}%
      \ifx \deflabelaux \deflabelbux \else \laberr{#1=(\deflabelaux)=(\deflabelbux)} \fi
      \fi
    \track{showlabel}{*}{#1}}

  \def \eqmark #1 {\advseqnumbering
    \eqno {(\current)}
    \deflabel{#1}{\current}}

  \def \subeqmark #1 {\global \advance\eqcntr by 1
    \edef\subeqmarkaux{\current.\number\eqcntr}
    \eqno {(\subeqmarkaux)}
    \deflabel{#1}{\subeqmarkaux}}

  \def \label #1 {\deflabel{#1}{\current}}
  \def \lcite #1{(#1\track{showlcit}{$\bullet$}{#1})}
  \def \forwardcite #1#2{\deflabel{#1}{#2}\lcite{#2}}


  \catcode`\@=11
  \def \c@itrk #1{{\bf #1}\track{showcitations}{\#}{#1}} 
  \def \c@ite #1{{\rm [\c@itrk{#1}]}}
  \def \sc@ite [#1]#2{[\c@itrk{#2}\hskip 0.7pt:\hskip 2pt #1]}
  \def \du@lcite {\if \pe@k [\expandafter \sc@ite \else \expandafter \c@ite \fi }
  \def \cite {\futurelet \pe@k \du@lcite }
  \catcode`\@=12


  \newcount \bibno \bibno = 0
  \def \newbib #1{\global\advance\bibno by 1 \edef #1{\number\bibno}}

  \def \bibitem #1#2#3#4{\smallbreak \item {[#1]} #2, ``#3'', #4.}

  \def \references {\begingroup \bigbreak \eightpoint
    \centerline {\tensc References}
    \nobreak \medskip \frenchspacing }


  \def \Headlines #1#2{\nopagenumbers
    \headline {\ifnum \pageno = 1 \hfil
    \else \ifodd \pageno \tensc \hfil \lcase {#1} \hfil \folio
    \else \tensc \folio \hfil \lcase {#2} \hfil
    \fi \fi }}

  \long \def \Quote #1\endQuote {\begingroup \leftskip 35pt \rightskip 35pt
\parindent 17pt \eightpoint #1\par \endgroup }
  \long \def \Abstract #1\endAbstract {\bigskip \Quote \noindent #1\endQuote}
  \def \Address #1#2{\bigskip {\tensc #1 \par \it E-mail address: \tt  #2}}
  
  \def \Note #1{\footnote {}{\eightpoint #1}}
  \def \Date #1 {\Note {\it Date: #1.}}


  \def \lcase #1{\edef \auxvar {\lowercase {#1}}\auxvar }

  \def \section #1 \par{\global \advance \secno by 1 \stno = 0
    \bigbreak \noindent {\bf \number \secno .\enspace #1.}
    \nobreak \medskip \noindent}

  \def \state #1 #2\par{\medbreak \noindent \advseqnumbering {\bf \current.\enspace #1.\enspace \sl #2\par }\medbreak }
  \def \definition #1\par {\state Definition \rm #1\par }

  \long \def \Proof #1\endProof {\medbreak \noindent {\it Proof.\enspace }#1
\ifmmode \eqno \endproofmarker $$ \else \hfill $\endproofmarker$ \looseness = -1 \fi \medbreak}

  \def \$#1{#1 $$$$ #1}
  \def\explain#1#2{\mathrel{\buildrel \hbox{\sixrm (#1)} \over #2}}
  \def \=#1{\explain{#1}{=}}

  \def \pilar #1{\vrule height #1 width 0pt}
  \def \stake #1{\vrule depth  #1 width 0pt}

  \newcount \footno \footno = 1
  \newcount \halffootno \footno = 1
  \def \footcntr {\global \advance \footno by 1
  \halffootno =\footno
  \divide \halffootno by 2
  $^{\number\halffootno}$}
  \def \fn#1{\footnote{\footcntr}{\eightpoint#1\par}}

  \begingroup
  \catcode `\@=11
  \global\def\eqmatrix#1{\null\,\vcenter{\normalbaselines\m@th\relax
      \ialign{\hfil$##$\hfil&&\kern 5pt \hfil$##$\hfil\crcr\relax
	\mathstrut\crcr\noalign{\kern-\baselineskip}\relax
	#1\crcr\mathstrut\crcr\noalign{\kern-\baselineskip}}}\,}
  \endgroup


  \def \Item #1{\smallskip \item {{\rm #1}}}
  \newcount \zitemno \zitemno = 0

  \def \izitem {\zitemno = 0}
  \def \zitemplus {\global \advance \zitemno by 1 \relax}
  \def \rzitem{\romannumeral \zitemno}
  \def \rzitemplus {\zitemplus \rzitem} 
  \def \zitem {\Item {{\rm(\rzitemplus)}}}
  \def \zitemmark #1 {\deflabel{#1}{\rzitem}}

  \newcount \nitemno \nitemno = 0
  
  \def \nitem {\global \advance \nitemno by 1 \Item {{\rm(\number\nitemno)}}}

  \newcount \aitemno \aitemno = -1
  \def \boxlet#1{\hbox to 6.5pt{\hfill #1\hfill}}
  \def \iaitem {\aitemno = -1}
  \def \aitemconv{\ifcase \aitemno a\or b\or c\or d\or e\or f\or g\or
h\or i\or j\or k\or l\or m\or n\or o\or p\or q\or r\or s\or t\or u\or
v\or w\or x\or y\or z\else zzz\fi}
  \def \aitem {\global \advance \aitemno by 1\Item {(\boxlet \aitemconv)}}
  \def \aitemmark #1 {\deflabel{#1}{\aitemconv}}


  \font\mf=cmex10
  \def\union {\mathop{\raise 9pt \hbox{\mf S}}\limits}
  \def\inters{\mathop{\raise 9pt \hbox{\mf T}}\limits}

  \def \<{\left \langle \vrule width 0pt depth 0pt height 8pt }
  \def \>{\right \rangle }
  \def \ds{\displaystyle}
  \def \and {\hbox {,\quad and \quad }}
  \def \calcat #1{\,{\vrule height8pt depth4pt}_{\,#1}}
  
  \def \imply {\mathrel{\Rightarrow}}
  \def \for #1{,\quad \forall\,#1}
  \def \endproofmarker {\square} 
  \def \"#1{{\it #1}\/}
  \def \inv {^{-1}}
  \def \*{\otimes}
  \def \caldef #1{\global \expandafter \edef \csname #1\endcsname {{\cal #1}}}
  \def \bfdef #1{\global \expandafter \edef \csname #1\endcsname {{\bf #1}}}
  \bfdef N \bfdef Z \bfdef C \bfdef R



  \catcode`\@=11

  \def\overparenOnefill{$\m@th
  \setbox0=\hbox{$\braceld$}%
  \braceld\leaders\vrule height\ht0 depth0pt\hfill
  \leaders\vrule height\ht0 depth0pt\hfill\bracerd$}

  \def\overparenOne#1{\mathop{\vbox{\m@th\ialign{##\crcr\noalign{\kern-1pt}
  \overparenOnefill\crcr\noalign{\kern3pt\nointerlineskip}
  $\hfil\displaystyle{#1}\hfil$\crcr}}}\limits}

  \def\overparenTwofill{$\m@th
  \lower 0.3pt \hbox{$\braceld$}
  \leaders\vrule depth 0pt height1pt \hfill
  \lower 0.3pt \hbox{$\bracerd$}$}

  \def\overparenTwo#1{\mathop{\vbox{\ialign{##\crcr\noalign{\kern-1pt}
  \overparenTwofill\crcr\noalign{\kern3pt\nointerlineskip}
  $\hfil\displaystyle{#1}\hfil$\crcr}}}\limits}

  \catcode`\@=12


  \def\X{X}
  \def\MyBlend{{\cal X}}

  \def\free{{\kern1pt*\kern1pt}}
  \def\genBlend{(A,B,i,j,\X)}
  \def\a{\alpha}
  \def\span#1{{\rm span}\ #1}
  \def\clspan#1{\overline{\rm span}\ #1}
  \def\id{\hbox{\sl id}}
  \def\E{E} \def\E{{\bf E}}
  \def\F{F} \def\F{{\bf F}}
  \def\Er{\E_*}
  \def\Fr{\F_*}
  \def\phir{\phi_*}

  \def\half{^{1/2}}
  \def\mhalf{^{-1/2}}
  \def\p{^\perp}   \def\p{^{^\perp}\kern-3pt}
  \def\op{^{op}}
  \font\rs = rsfs10 \def\L{\hbox{\rs L}\, }
  \def\vezes{\circledast}
  \def\labelarrowlen#1#2{\ {\buildrel #1\over {\hbox to #2{\rightarrowfill}}}\ }

  \def\vec#1{\left(#1\pilar{9pt}\right)_n}
  \def\f{^\infty}%
  \def\tilde{\widetilde}
  \def\array#1{\left[\matrix{#1}\right]}
  \def\m#1#2#3#4{\array{#1&#2\cr#3&#4}}
  \def\d#1#2{\m{#1}{}{}{#2}}
  \def\tensor#1{\!\!\mathop{\*}\limits_{#1}\!\!}
  \def\tmin{\tensor{min}}
  \def\tmax{\tensor{max}}
  \def\soma{\sum_{i=1}^n}
  \def\subsetequp{\cup \kern 1pt \vrule height 6pt depth 0pt width 0.5pt}
  \def\lt{L^2(S^1)}
  \def\cs{C(S^1)}

  \def\mybib#1#2{\global \expandafter \edef \csname #1\endcsname {#2}}
  \def\mybib#1#2{\global \advance \bibno by 1 \global \expandafter \edef \csname #1\endcsname {\number\bibno}}


  \mybib{AVRR}         {AVRR}
  \mybib{BCM}          {BCM}
  \mybib{Brz}          {Brz}
  \mybib{DT}           {DT}
  \mybib{ED}           {ED}
  \mybib{EN}           {EN}
  \mybib{Gucc}         {Gucc}
  \mybib{Ped}          {Ped}
  \mybib{Oka}          {O}
  \mybib{PP}           {PP}
  \mybib{Po}           {Po}
  \mybib{Rieffel}            {R1}
  \mybib{RifRot}            {R2}
  \mybib{PV}            {PV}
  \mybib{W}            {W}

  \def\tittext{BLENDS AND ALLOYS}
  \Headlines {\tittext}{R.~Exel}
  \null\vskip -1cm
  \centerline{\titlefont \tittext}

  \bigskip
  \centerline{\tensc Ruy Exel}

  \Date{25 March 2012}

  \footnote{\null}
  {\eightrm 2010 \eightsl Mathematics Subject Classification:
  \eightrm  Primary
	46L05, 
	16S99; 
    Secondary
	16S35, 
	16S40. 
  }

  \Note
  {\it Key words and phrases: \rm Algebra, C*-algebra, crossed
product, tensor product, algebra structure, conditional expectation,
index finite type, Jones' basic construction, commuting square, blend,
alloy.}

  \bigskip

  \Abstract
  Given two algebras $A$ and $B$, sometimes assumed to be C*-algebras,
we consider the question of putting algebra or C*-algebra structures on the tensor
product $A\*B$.  In the C*-case,  assuming $B$ to be two-dimensonal,
we characterize all possible such C*-algebra structures in terms of an
action of the cyclic group $\Z_2$.  An example related to
commuting squares is also discussed.
  \endAbstract

  \section Introduction

When $G$ is a group and
  $
  \alpha : G \to {\rm Aut}(A)
  $
  is an action of $G$ on a unital $K$-algebra $A$, one may form the
the \"{crossed product algebra} (also known  among algebraists as the \"{skew group
algebra})\def\cp{A\rtimes_\alpha G}
  $
  \cp.
  $
  As a vector space $\cp$ is just the tensor product $A\*K(G)$,
where $K(G)$ denotes the group algebra of $G$ with coefficients in the
base field $K$.  The multiplication
operation on $\cp$ is given by
  $$
  (a\*g)  (b\*h) = a\alpha_g(b)\*gh
  \for a, b\in A \for g, h\in G.
  $$

Researchers working with crossed products are used to thinking that
the above multiplication operation on $A\*K(G)$ has been \"{twisted},
or \"{skewed} by the group action in relation to the usual tensor
product multiplication.  Viewing things from this point of view, one
can't help but to ask in how many other ways can the usual
mutiplication on a tensor product algebra
  be similarly modified.

Returning to the example of crossed products above, it turns out that the maps
  $$
  i:a\in A\mapsto a\*1\in \cp
  $$and$$
  j:y\in K(G)\mapsto 1\*y\in \cp
  $$
  are algebra homomorphisms and, in addition, $\cp$ is equal to the linear span of both
  $
  i(A)j\big(K(G)\big) $ and $ j\big(K(G)\big)i(A).
  $

 Motivated by these properties we define a \"{blend} of algebras as
being a quintuple
  $$
  \MyBlend =\genBlend,
  $$
  where $A$, $B$ and $\X$ are unital algebras (see below for the
definition in the non-unital case), and
  $$
  i: A \to \X
  \and
  j: B \to \X
  $$
  are unital homomorphisms such that the maps
  $$
  i\vezes j\ : \ a\*b \in A\*B \quad \longmapsto \quad i(a)j(b) \in \X,
  $$
  and
  $$
  j\vezes i\ : \ b\*a \in B\*A \quad \longmapsto \quad j(b) i(a)\in \X
  $$
  are surjective.  If moreover
$i\vezes j$ and   $j\vezes i$ are one-to-one, we say that $\MyBlend$ is an
\"{alloy}.

Working in the category of C*-algebras we introduce similar notions,
but the above requirement that $i\vezes j$ and $j\vezes i$ are
surjective is replaced by the weaker requirement that they have dense range.

  Given an algebraic alloy, we may identify $A\*B$ with $\X$ under
$i\vezes j$, and hence make $A\*B$ an algebra.  Its multiplication
operation will therefore satisfy
  $$
  (a'\*1)(a\*b) = a'a\*b
  \and
  (a\*b)(1\*b') = a\*bb'.
  \eqmark BasicRules
  $$
  This is evidently  not enough to characterize the whole multiplication
operation in $A\*B$,  but if one is also given the map
  $$
  \tau : B\*A  \to A\*B,
  $$
  defined by $\tau = (i\vezes j)\inv(j\vezes i)$,
  then the product between $1\*b$ and $a\*1$ may be written in terms
of $\tau$ as
  $$
  (1\*b)(a\*1) = \tau(b\*a).
  $$
  In fact the multiplication operation on
$A\*B$ may be completely recovered,  given $\tau$, as follows:
  $$
  \matrix{ (a_1\*b_1) (a_2\*b_2)
  & = & (a_1\*1) (1\*b_1) (a_2\*1) (1\*b_2) \cr
  \pilar{12pt}  & = & (a_1\*1)\hfill \tau(b_1\*a_2)\hfill (1\*b_2)}
  \eqmark TheProduct
  $$
  while the final answer may be reached upon an application of \lcite{\BasicRules}.

  Conversely, given a map $\tau : B\*A  \to A\*B$, one may take the above
answer as the \"{definition} of a multiplication operation on $A\*B$,
which will evidently not always be associative but,
at least in theory, one may spell out a condition on $\tau$ for the associativity
to hold.

At this point I must confess that I could not find any nice looking
associativity condition.  This is perhaps an indication that this
whole circle of ideas is thornier than one would initially believe.
For starters, as already seen, this construction would encompass all
group actions on algebras!

Of course  similar thorny  questions may be also asked in case $A$ and $B$ are C*-algebras,
but then one wouldn't expect to get a C*-algebra right away, unless $A\*B$ is
completed under some suitable norm.

In the purely algebraic situation this question has already been
extensively treated, but almost always $B$ is supposed to have some
extra structure, like that of a Hopf algebra \cite {\BCM}, \cite {\DT} or some
similar structure \cite{\Gucc}, \cite{\AVRR}.

Perhaps the most general approach to this problem is due to
Brzezi\'nski \cite{\Brz}, where no algebra structure whatsoever is assumed
on $B$,  which is only assumed to be a vector space!

Among other things the beauty of Mathematics rests on the fact that, no matter
how hard is the problem facing us, there is always an easier special
case which may be effectively studied and which, one hopes, may lead the way to
new ground.

Far from attempting a complete theory of blends or alloys in either
the purely algebraic or the C*-algebraic context, the aim of this
paper is to exploit a few concrete situations in which we have found
some surprising amount of mathematical structure, but  which apparently have
not yet been discussed in the literature.

In our first example, a truly humbling experience, we study all
possible C*-algebra structures on $A\*B$, when $A$ is an arbitrary
unital C*-algebra and $B=\C^2$.  Curiously this turns out to be quite
involving and the answer is that crossed products of $A$ by actions of
$\Z_2$ provide all possible examples!

This evidently begs for a generalization to $B=\C^n$, for $n>2$, but
unfortunately our methods do not seem to extend beyond the case $n=2$.

Our second main example is based on Jones' basic construction, as
generalized by Watatani \cite{\W} to C*-algebras.  As initial data we
consider a commuting square \cite{\Po} of C*-algebras, meaning
C*-algebras $A$, $B$, $C$ and $D$, such that
  $$
  \matrix{
  A & \supseteq & B \cr
  \vrule height 14pt depth 10pt width 0pt \subsetequp   &&   \subsetequp  \cr
  C & \supseteq & D,
  }
  $$
  and $D=B\cap C$.  We also suppose we are given conditional expectations
  $$
  E:A\to B \and F:A\to C,
  $$
  satisfying $EF=FE$.   Evidently one then has that
  $
  G := E F
  $
  is a conditional expectation onto $D$.
  Given such a commuting square we let $M$ be the Hilbert $D$-module
obtained by completing $A$ under the $D$-valued inner-product
  $$
  \langle a,b\rangle = G(a^*b), \quad \forall a,b\in A.
  $$
  We moreover let $e$, $f$ and $g$ be the projections on $M$
obtained by extending $E$, $F$ and $G$, respectively, to $M$.
Viewing the left action of $A$ on $M$ as a homomorphism $\lambda: A
\to \L(M)$,  where $\L(M)$ denotes the algebra of all adjointable  operators on $M$, we
introduce our main players, namely the C*-algebras
  $$
  K_g = \clspan{\lambda(A)g\lambda(A)}, \quad
  K_e = \clspan{\lambda(A)e\lambda(A)} \and
  K_f = \clspan{\lambda(A)f\lambda(A)}.
  $$
  We then show that there are natural maps
  $$
  i_e:K_e \to M(K_g) \and
  i_f:K_f \to M(K_g),
  $$
  where $M(K_g)$ stands for the multiplier algebra of $K_g$.
  The main oustanding question left unresolved by this work,
incidentally also the question that motivated me to consider the
notion of a C*-blend, is whether or not
  $$
  \big(K_e,K_f,i_e,i_f,K_g\big)
  $$
  is a C*-blend (please see below for the definition of a C*-blend in
the non-unital case).

  Having little to say about this problem in its full generality we
consider two special cases where we are able to give positive
answers.  The first one is the case in which $G$  is of index finite type
\cite{\W},  while in the second case we assume that $A$ is a unital
commutative algebra and $D = \C\, 1$.

After having discussed a preliminary notion of blends and alloys with George
Elliott and Zhuang Niu, they have found an interesting example of
a C*-alloy related to the irrational rotation C*-algebra, which we
briefly describe below with their kind permission.

Last but not least I would like to thank V.~Jones and G.~Elliott for stimulating
conversations  while this work was being prepared.

\section General concepts

In this section we shall describe the main concepts to be treated in
this work and which will be further developped in several different
directions in the forthcoming sections.

We will assume througout that $R$ is a commutative unital
ring, quite often specialized to be the field of complex numbers.  By
an algebra we will always understand an associative $R$-algebra.

If $A$ is an $R$-algebra we will denote by $M(\X)$ the algebra of
multipliers of $\X$ (see e.g. \cite[2.1]{\ED} for a definition).  In badly behaved cases
$\X$ is not necessarily isomorphic to a subalgebra of $M(\X)$, so we
will assume that the natural map from $\X$ to $M(\X)$ is injective and
hence we will think of $\X$ as a subalgebra of $M(\X)$.  Please see the
discussion following \cite[2.2]{\ED}.

\definition  \label AlgDefi Consider a quintuple
$\MyBlend = \genBlend$, where $A$,  $B$ and $\X$ are $R$-algebras, and
  $$
  i: A \to M(\X)
  \and
  j: B \to M(\X)
  $$
  are homomorphisms.  Also consider the $R$-linear
maps
  $$
  i\vezes j\ : \ a\*b \in A\*_RB \quad \longmapsto \quad i(a)j(b) \in M(\X), $$  and  $$
  j\vezes i\ : \ b\*a \in B\*_RA \quad \longmapsto \quad j(b) i(a)\in M(\X).
  $$
  We will say that $\MyBlend$ is:
  \iaitem
  \aitem a \"{blend} when the ranges of both $i\vezes j$
and $j\vezes i$ coincide with $\X$ (or rather its cannonical copy
within $M(\X)$),
  \aitem an \"{alloy} when, in addition to (a),  one has
that
  $i\vezes j$ and $j\vezes i$ are one-to-one.
  \smallskip \noindent
  In either case, when $A$, $B$, and $\X$ are unital algebras, and both
$i$ and $j$ preserve the unit, we will say that $\MyBlend$ is a \"{unital}
blend or alloy.

As already mentioned in the introduction, if $G$ is a group
acting on a unital algebra $A$, then
  $$
  \big(A,R(G),i,j,A\rtimes G\big)
  $$
  is an example of a blend,  which is easily seen to be an alloy.

\definition Let
  $\MyBlend_1 = (A,B,i_1,j_1,\X_1)$ and $\MyBlend_2 = (A,B,i_2,j_2,\X_2)$ be
blends.  A \"{morphism} from $\MyBlend_1$ to $\MyBlend_2$ is a homomorphism
  $$
  \phi:M(\X_1) \to M(\X_2),
  $$
  such that
  $\phi i_1= i_2$, and $\phi j_1= j_2$.  The morphism will
be termed an \"{isomorphism} when $\phi$ is bijective.  If $\MyBlend_1$ and
$\MyBlend_2$ are unital blends, we will say that $\phi$ is a \"{unital
morphism} provided $\phi(1) = 1$.

Given a morphism $\phi$, as above, notice that
  $$
  \phi(i_1\vezes j_1) = i_2\vezes j_2.
  \eqmark MorphismAndVezes
  $$
  Therefore it is evident that  $\phi(\X_1) = \X_2$.

\state Proposition \label BlendToAlloy Every morphism from a blend to
an alloy is an isomorphism.

\Proof
  Let $\MyBlend_1 = (A,B,i_1,j_1,\X_1)$ be a blend,  let $\MyBlend_2 =  (A,B,i_2,j_2,\X_2)$
be an alloy,  and let $\phi$ be a morphism from $\MyBlend_1$ to $\MyBlend_2$.
Denote by $\psi$ the restriction of $\phi$ to $\X_1$, viewed as a
map onto $\X_2$,  and observe that by \lcite{\MorphismAndVezes} we have that
  $$
  \psi(i_1\vezes j_1) = i_2\vezes j_2.
  $$
  By hypothesis we have that $i_2\vezes j_2$ is one-to-one, and hence
so is $i_1\vezes j_1$.   Thus both $i_1\vezes j_1$ and $i_2\vezes j_2$
are bijective maps,  respectively onto $\X_1$ and $\X_2$.
Therefore
  $$
  \psi= (i_2\vezes j_2)(i_1\vezes j_1) \inv,
  $$
  and hence $\psi$ is bijective, and the reader may now prove that
$\phi$ is necessarily bijective as well.
  \endProof 

\section C*-blends

Let us now consider C*-algebraic versions of the above concepts, and
hence we now assume that our base ring $R$ is the field of complex numbers.

\definition Consider a quintuple
$\MyBlend = \genBlend$, where $A$,  $B$ and $\X$ are C*-algebras, and
  $$
  i: A \to M(\X)
  \and
  j: B \to M(\X)
  $$
  are *-homomorphisms.  Also consider the linear
maps
  $i\vezes j$ and   $j\vezes i$ described in \lcite{\AlgDefi}
  (where $A\*_{\C} B$ and $B\*_{\C}A$ refer to the algebraic tensor product over
the field of complex numbers).
  We will say that $\MyBlend$ is:
  \iaitem \aitem a \"{C*-blend} when the ranges of both $i\vezes j$
and $j\vezes i$ are contained and dense in $\X$,
  \aitem a \"{C*-alloy} when, in addition to (a),  one has
that
  $i\vezes j$ and $j\vezes i$ are one-to-one.
  \smallskip \noindent
  In either case, when $A$, $B$, and $\X$ are unital algebras, and both
$i$ and $j$ preserve the unit, we will say that $\MyBlend$ is a \"{unital}
C*-blend or C*-alloy.

Observe that, since the range of $\,j\vezes i\,$ is the adjoint of the
range of $\,i\vezes j\,$, if the range of one of this maps is
contained and dense in $\X$, then so is the other.

The definition of a \"{C*-alloy} given above should be considered as
tentative for the following reason: given C*-algebras $A$ and $B$,
think of the case of the minimal versus the maximal tensor products,
here denoted $A\tmin B$ and $A\tmax B$, respectively.
  One evidently has two natural C*-alloys, namely
  $$
  \MyBlend_{min} = (A,B,i,j,A\tmin B) \and   \MyBlend_{max}=(A,B,i,j,A\tmax B).
  $$
  However, as an alloy, $\MyBlend_{min}$ does not satisfy any sensible
generalization of \lcite{\BlendToAlloy}, since the natural map
  $$
  \phi:A\tmax B\to A\tmin B
  $$
  provides a morphism of C*-blends which is not always an isomorphism.
On the other hand it would be interesting to check if $\MyBlend_{max}$
satisfies anything like \lcite{\BlendToAlloy}.

Note that a C*-blend is not necessarily an algebraic blend, since
the requirement of having dense range is weaker than that of being
surjective.  However, occasionally we shall be interested in a concept
which subsumes both.

\definition A C*-blend $\genBlend$ will be called \"{strict} if it also
satisfies the conditions of \lcite{\AlgDefi.a}, that is, if $i\vezes j$
and $j\vezes i$ are onto $\X$.

In all of our uses of the notion of strict C*-blends below, one of the
algebras involved will be finite dimensional.

\state Remark \label ShortNotationBlend \rm In the unital case, if
$\genBlend$ is either an algebraic or C*-alloy, we may identify $A$
and $B$ with the subalgebras of $\X$ given by $i(A)$ and $j(B)$,
respectively.
  When such an identification is being made we will sometimes use the
short hand notation $(A,B,\X)$ to indicate this.

An interesting source of examples of C*-blends is obtained from the
theory of crossed products.  Given a locally compact group $G$ acting
on a C*-algebra $A$, consider the full
crossed product\fn{Similar conclusions may be obtained if we choose
the reduced crossed product instead.}
  $$
  A\rtimes G.
  $$
  Representing   $A\rtimes G$ faithfully on a Hilbert space $H$,  one
deduces from \cite[7.6.4]{\Ped} that there exists a covariant representation
  $(\pi,u)$, such that the given faithfull representation coincides
with $\pi\times u$.

  It is well known \cite[7.6.6]{\Ped} that the range of $\pi$ lies in
the multiplier algebra of $A\rtimes G$.  Moreover, if
  $$
  \rho : C^*(G) \to B(H)
  $$
  denotes the integrated form of $u$, then the range of $\rho$ is also
contained in $M(A\rtimes G)$.

\state Proposition \label CrossIsBlend
  $
  \big(A,C^*(G),\pi,\rho,A\rtimes G\big)
  $
  is a C*-blend.

\Proof Let us first prove that the range of $\pi\vezes\rho$ is
contained in $A\rtimes G$.  For this let $a\in A$ and $f\in C^*(G)$,
and write $f = \lim_nf_n$, with $f_n\in C_c(G)$.

  Viewing $a\*f_n$ as an element of $C_c(G,A)$, and observing
that
  $$
  (\pi\times u)(a\*f_n) = \pi(a)\rho(f_n),
  $$
  by \cite[7.6.4]{\Ped}, one sees that $\pi(a)\rho(f_n)$ is in
$A\rtimes G$.
So
  $$
  (\pi\vezes\rho)(a\*f) = \pi(a)\rho(f) =
  \lim_n\pi(a)\rho(f_n) \in A\rtimes G,
  $$
  as claimed.

  Next, let us check that the range of $\pi\vezes\rho$ is dense in
$A\rtimes G$.  For this, let us view $A\*C_c(G)$ as a subspace of
$C_c(G,A)$ in the usual fashion.  The former may be easily shown to be
dense in the latter under the inductive limit topology.

On the other hand, by the very definition of crossed products, the
image of $C_c(G,A)$ under $\pi\times u$ is a dense subalgebra of
$A\rtimes G$.
  Therefore
  $$
  (\pi\vezes\rho)\big(A\*C_c(G)\big)  =
  \pi(A) \rho\big(C_c(G)\big) = (\pi\times u)\big(A\*C_c(G)\big)
  $$
  is dense in $A\rtimes G$.  Evidently we also have that
  $
  (\pi\vezes\rho)\big(A\*C^*(G)\big)
  $
  is dense.

  Upon taking adjoints, the same conclusions apply to
  $\rho\vezes\pi$, concluding the proof.
  \endProof

  If $G$ is a finite group one may in fact prove the above to be a
strict C*-alloy.  In this case, one has that $C^*(G)$ is finite
dimensional and this might lead one to expect C*-blends to be strict
whenever one of the algebras involved is finite dimensional.  However,
in section \forwardcite{BadExample}{9} we shall present an example of
a non strict C*-blend $\genBlend$, in which $B$ is a two-dimensional
algebra.

\section The Elliott-Niu alloy

In this short section we describe,  without proofs, an interesting example  of C*-alloy
found by Elliott and Niu.

Given an irrational number $\theta$, consider the unitary operators
$u$ and $v$ on $\lt$ defined by
  $$
  u(\xi)\calcat z = z\xi(z)
  \and
  v(\xi)\calcat z = \xi(e^{-2\pi i\theta z}),
  $$
  for all $\xi\in \lt$ and $z\in S^1$.  It is well known that the
C*-algebra generated by $u$
and $v$ is isomorphic to the irrational rotation C*-algebra
$A_\theta$.  Since the latter is a crossed product of $\cs$ by
$\Z$,  we may use \lcite{\CrossIsBlend} to obtain a blend (cf.~Remark \lcite{\ShortNotationBlend})
  $$
  \big(\cs,\cs, A_\theta\big),
  $$
  where the second occurence of $\cs$ above is justified by the fact that
it is isomorphic to $C^*(\Z)$.
  Let $\chi$ be the characteristic function of the arc
  $$
  J = \{e^{it}: 0 \leq t \leq 2\pi\theta\},
  $$
  and let $p$ be the operator on $\lt$ given by
  $$
  p(\xi)\calcat z = \chi(z)\xi(z),
  $$
  so that $p$ is just the spectral projection of $u$ associated to the
arc $J$.

  Denote by $B$ the closed *-algebra of operators on $\lt$ generated
by the set
  $$
  \{u,v^kpv^{-k}: k\in \Z\}.
  $$
  It is not hard to show that $B$ is isomorphic to the algebra
  $$
  C(S^1;\theta),
  $$
  formed by all bounded, complex valued functions on the circle which are
continuous at all points except possibly for points in the orbit of 1
under rotation by $2\pi\theta$,  where the lateral limits exist.

Since $B$ is clearly invariant under rotation by $2\pi\theta$
(conjugation by $v$), one may also form the crossed product
$B\rtimes\Z$ and, as above, it is not surprising that
  $$
  \big(B,\cs, B\rtimes\Z\big)
  $$
  is also a blend.  The interesting aspect of the Elliott-Niu example
is that it extrapolates the realm of crossed products as follows.

Let
$q$ be the spectral projection of $v$ (as opposed to $u$) associated to
the arc   $J$, and let $C$ be the
closed *-algebra of operators on $\lt$ generated
by the set
  $$
  \{v,u^kqu^{-k}: k\in \Z\}.
  $$
  Again it is possible to prove that $C$ is isomorphic to $C(S^1;\theta)$.

\state Theorem {\rm (Elliott-Niu \cite{\EN})} The set
  $$
  X_\theta = \clspan \{bc: b\in B,\ c\in C\}
  $$
  is a subalgebra of operators on  $\lt$, and consequently
  $(B,C,X_\theta)$ is a C*-blend, which is in fact a C*-alloy.

A proof of this result will hopefully appear soon.

In analogy with the highly influential isomorphism problem for irrational rotation
C*-algebras (\cite{\RifRot}, \cite{\PV}), one may ask:

\state Question \rm Given real numbers $\theta$ and $\theta'$,  when
are $X_\theta$ and $X_{\theta'}$ isomorphic?

\section Algebraic  alloys with $B=R^2$

The examples based on crossed products above, together with the
complexity and wide variety of group actions on C*-algebras, should be
enough to convince the reader that the concept of a blend, in both the
algebraic and the C* version, is a rather general and deep one.
Understanding it in full is therefore a colossal task not likely to be
accomplished in the near future.

In order to be able to say anything meaninful about blends and alloys we will
make a drastic simplification by assuming that one of the algebras
involved is as elementary as can possibly be, namely $B=\C^2$.  In
this case we will give a complete classification of unital strict C*-alloys.

In preparation for this we shall now study the algebraic counterpart,
namely alloys in which $B$ is the direct sum of two copies of the
coefficient ring $R$.

As the above title suggests we shall fix, for the remainder of this
section, a unital alloy $\genBlend$ such that $B=R^2$, the latter viewed as an algebra with the usual
coordinatewise operations.
As already mentioned we  will think of $A$ and $R^2$ as subalgebras of
$\X$ and we will refer to our alloy simply as $(A,R^2,\X)$.
We will denote by $p$ and $q$ the standard idempotents of $R^2$, namely
  $$
  p = (1,0) \and q = (0,1).
  $$

\state Proposition \label UniqueDecomp Every element $c\in \X$ admits
unique decompositions as
  $$
  c = ap+bq,
  \subeqmark UniqueLeft
  $$
  and
  $$
  c = pa'+qb',
  \subeqmark UniqueRight
  $$
  where $a, b, a', b'\in A$.

\Proof This follows from the corresponding statements about decomposing an
element in $A\*_RR^2$ as $a\*p + b\*q$, and similarly with respect to
  $R^2\*_RA$.
  \endProof

The structure of the multiplication operation on $\X$ is encoded by two crucial operators
on $A$ introduced in the following:

\state Proposition \label EnterEnadF There are unique linear maps $\E,\F:A\to A$,
such that
  $$
  pap = \E(a)p
  \and
  qaq = \F(a)q,
  $$
  for all $a\in A$.
  In addition, denoting by $\E\p=\id_A-\E$, and $\F\p=\id_A-\F$, we have that
  \izitem
  \zitem $pa = \E(a)p + \F\p(a)q$,
  \zitem $qa = \E\p(a)p + \F(a)q,$
  \medskip\noindent
  for all $a\in A$.

\Proof
  Given $a\in A$, we may use \lcite{\UniqueLeft} to write
  $$
  pa = x_ap+y_aq,
  \subeqmark PxaPbQ
  $$
  for a unique pair $(x_a, y_a)\in A^2$.  The map
  $$
  \E:a\in A\mapsto x_a\in A
  $$
  is therefore well defined and may be easily proven to be linear. Right
multiplying \lcite{\PxaPbQ} by $p$ then leads to
  $$
  pap = x_ap = \E(a)p.
  $$
  The uniqueness of $\E$ follows immediately from \lcite{\UniqueDecomp} and a
similar reasoning applies to prove the existence and uniqueness of $\F$.
  As for the last part of the statement, given $a$ in $A$, write
  $$
  pa = x_ap + y_aq \and qa = z_ap + w_aq.
  $$
  Clearly $x_a = \E(a)$, and $w_a = \F(a)$, while
  $$
  ap + aq = a (p+q) =
  a = (p + q)a = pa + qa = \big(\E(a) + z_a\big)p + \big(y_a+ \F(a)\big)q.
  $$
  Again by the uniqueness part of \lcite{\UniqueDecomp} we deduce that
  $$
  z_a = a-\E(a) = \E\p(a)
  \and y_a = a- \F(a) = \F\p(a).
  $$
  This concludes the proof.
  \endProof

  The fact that $\E$ and $\F$ really do encode the multiplication operation on $\X$
is made clear by the following:

  \state Proposition \label FormulaForMult Given $c_1=a_1p+b_1q$, and
$c_2=a_2p+b_2q$ in $\X$, one has that
  $$
  c_1c_2 =
  \big(a_1\E(a_2) + b_1\E\p(a_2)\big)p + \big(a_1\F\p(b_2)+b_1\F(b_2)\big)q.
  $$

  \Proof Left for the reader. \endProof

Because of the importance of these maps in describing the algebra structure of $\X$ they ought
to be given a name.

\definition  We will refer to $(\E,\F)$ as the \"{left intrinsic pair} for the alloy $(A,R^2,\X)$.

  Some further important properties of $\E$ and $\F$ are studied next.

\state Proposition \label ManyRels Both $\E$ and $\F$ are idempotent operators and moreover, for
every $a$ and $b$ in $A$, one has that
  \izitem
  \zitem $\E(ab) = \E(a)\E(b) + \F\p(a)\E\p(b)$
  \zitem $\F(ab) = \E\p(a)\F\p(b) + \F(a)\F(b)$
  \zitem $\E\p(ab) = \E\p(a)\E(b) + \F(a)\E\p(b)$
  \zitem $\F\p(ab) = \E(a)\F\p(b) + \F\p(a)\F(b)$

\Proof
  Given $a$ in $A$, one has that
  $$
  \E(a) p = pap = p(pap) = p\E(a)p = \E(\E(a))p,
  $$
  whence $\E^2=\E$, and one similarly proves that $\F^2=\F$.
  In order to prove (i) notice that
  $$
  \E(ab)p = (pa)bp =
  \big(\E(a)p + \F\p(a)q\big)bp =
  \E(a)\E(b)p + \F\p(a)\E\p(b)p \$=
  \big(\E(a)\E(b) + \F\p(a)\E\p(b)\big)p,
  $$
  proving (i).  The proof of (ii) follows similar lines.  In order to prove
(iii) we start with its right hand side:
  $$
  \E\p(a)\E(b) + \F(a)\E\p(b) =
  \big(a-\E(a)\big)\E(b) + \big(a-\F\p(a)\big)\E\p(b) \$=
  a\E(b) -\E(a)\E(b) + a\E\p(b) -\F\p(a)\E\p(b) \$=
  ab -\E(a)\E(b) -\F\p(a)\E\p(b) \={i}
  ab - \E(ab) = \E\p(ab).
  $$
  One may similarly show that (ii) implies (iv).
  \endProof

Another crucial piece of information to be extracted from
$(A,R^2,\X)$ is as follows:

\state Proposition \label EnterPhi The map
  $$
  \phi:  a\in A \mapsto  \E(a) + \F(a) - a\in A
  $$
  is multiplicative.

\Proof
  Given $a$ and $b$ in $A$ we have
  $$
  \phi(ab) = \E(ab) - (ab - \F(ab)) = \E(ab)-\F\p(ab) \={\ManyRels.i \& iv} $$$$ =
  \E(a)\E(b) + \F\p(a)\E\p(b)  -\E(a)\F\p(b) - \F\p(a)\F(b) \$=
  \E(a)\big(\E(b)  -\F\p(b)\big)+ \F\p(a)\big(\E\p(b)  - \F(b)\big) =
  \E(a)\phi(b) - \F\p(a)\phi(b) \$=
  \big(\E(a)- \F\p(a)\big)\phi(b) = \phi(a)\phi(b).
  \endProof

It is our next immediate goal to show that $\phi$ is indeed an automorphism of
$A$.  In order to do so we need to introduce the right-handed versions of the operators
$\E$ and $\F$.

\state Proposition \label EnterDualEandF There are unique linear maps $\Er,\Fr:A\to A$, such that
  $$
  pap = p\Er(a)
  \and
  qaq = q\Fr(a)
  \for a\in A.
  $$
  In addition,
  \izitem
  \zitem $ap = p\Er(a) + q\Fr\p(a)$,
  \zitem $aq = p\Er\p(a) + q\Fr(a),$
  \medskip\noindent  for all $a\in A$.

\Proof By considering the oposite algebras we may think of the alloy $(A\op,R^2,\X\op)$ (since
$R^2$ is commutative, it coincides with its oposite algebra).  We may then apply
\lcite{\EnterEnadF} to obtain the corresponding versions of $\E$ and $\F$, which we denote by
$\Er$ and $\Fr$, respectively.  The conditions in the statement are then obtained from the
corresponding conditions in \lcite{\EnterEnadF}, once the order of the factors in all products
are suitably reversed.  \endProof

\definition We shall refer to $(\Er,\Fr)$ as the \"{right intrinsic
pair} for the alloy $(A,R^2,\X)$.

The following result lists some relevant relations satisfied by $\E$ and $\F$ together with their
right-handed versions.

\state Proposition \label JointRel The maps $\E$, $\F$, $\Er$ and $\Fr$ introduced  above satisfy:
  \izitem
  \zitem $\E\Er = \E$,
  \zitem $\F\Er = \Er$,
  \zitem $\F\Fr = \F$,
  \zitem $\E\Fr = \Fr$.

\Proof
  Given $a$ in $A$ we have that
  $$
  \E(a)p = pap =
  p\Er(a) =
  \E\big(\Er(a)\big)p + \F\p\big(\Er(a)\big)q.
  $$
  By  \lcite{\UniqueDecomp} we then deduce that
  $\E = \E\Er$, and
  $$
  0 = \F\p\big(\Er(a)\big) = \Er(a) - \F\big(\Er(a)\big),
  $$
  whence $\Er = \F\Er$.  The proofs of (iii) and (iv) are similar.
  \endProof

\state Proposition \label InverseOfPhi The map $\phi$ introduced in \lcite{\EnterPhi} is an
automorphism of $A$  and its inverse is given by
  $$
  \phir := \Er+\Fr-\id_A.
  $$

\Proof
  We have that
  $$
  \phi\phir = (\E+\F-\id_A)(\Er+\Fr-\id_A) \$=
  \E\Er+\E\Fr-\E + \F\Er+\F\Fr-\F - \Er-\Fr+\id_A \={\JointRel} $$$$=
  \E  + \Fr-\E +  \Er+\F  -\F - \Er-\Fr+\id_A  = \id_A,
  $$
  and one may similarly prove that   $\phir\phi = \id_A$.
\endProof

\definition The automorphism $\phi$ above shall be called the \"{intrinsic
automorphism} of the alloy $(A,R^2,\X)$.

It is interesting that $\phi$ may be extended to the whole of $\X$, as shown in the next:

\state Proposition \label EnterPsi There exists a unique automorphism $\Phi$ of $\X$, extending
$\phi$, and such that $\Phi(p) = q$, and $\Phi(q) = p$.

\Proof  Given $c$ in $\X$,  let $c=ap+bq$ be its unique decomposition according
to \lcite{\UniqueDecomp} and put
  $$
  \Phi(c) = \phi(b)p + \phi(a)q.
  $$
  In order to show that $\Phi$ is multiplicative we first notice that
  $$
  \phi \E = \F\phi\and \phi \F = \E\phi,
  \subeqmark PhiEEqFFhi
  $$
  as the reader may easily verify by writing $\phi = \F-\E\p = \E-\F\p$.
  Given $c_1$ and $c_2$ in $A$,
write
  $c_i=a_ip+b_iq$, for $i=1,2$, so that
  $$
  \Phi(c_1c_2) \={\FormulaForMult}
  \phi\big(b_1\F(b_2)+a_1\F\p(b_2)\big)p + \phi \big(a_1\E(a_2) + b_1\E\p(a_2)\big) q \$=
  \Big(\phi(b_1)\E\big(\phi(b_2)\big) +\phi(a_1)\E\p\big(\phi(b_2)\big)\Big)p +
  \Big(\phi(a_1) \F\big(\phi(a_2)\big) + \phi(b_1)\F\p\big(\phi(a_2)\big)\Big)q \={\FormulaForMult} $$$$=
  \big(\phi(b_1)p+\phi(a_1)q\big)   \big(\phi(b_2)p+\phi(a_2)q\big) =
  \Phi(c_1)  \Phi(c_2).
  $$
  This shows that $\Phi$ is an endomorphism of $\X$.

  Glancing at the definition of $\Phi$ one immediately checks that it is injective and surjective
and hence that $\Phi$ is indeed an automorphism.  Uniqueness of $\Phi$ is also evident.
  \endProof

\section Conditional Expectations

We continue enforcing the standing hypothesis set out at the beginning
of the previous section, namely that $(A,R^2,\X)$ is a fixed unital algebraic
alloy.

From now on we shall be interested in a concept borrowed from probability theory:

\definition \label ConExp By a \"{conditional expectation} from $\X$ to $A$ we shall mean a linear map
  $$
  G: \X\to A,
  $$
  which is an $A$-bimodule map and which coincides with the identity on $A$.

\state Proposition \label HvsCondExp There is a one-to-one correspondence between the set of all
conditional expectations $G$ from $\X$ to $A$ and the set of all elements $h$ in $A$
such that
  $$
  ha =   \phi(a)h + \F\p(a)
  \for a \in A.
  \subeqmark haphihfp
  $$
  The correspondence is given by $G\mapsto h = G(p)$.

\Proof
  We first claim that, for every $a\in A$,
we have that
  $$
  pa =  \phi(a)p + \F\p(a).
  \subeqmark paphiFp
  $$
  This follows from
  $$
  pa = \E(a)p+\F\p(a)(1-p) =
  \big(\E(a)-\F\p(a)\big)p + \F\p(a) =
  \phi(a)p + \F\p(a).
  $$
  Given a conditional expectation $G$, let $h = G(p)$.  Then
  $$
  ha=G(p)a = G(pa) \={\paphiFp}
  G\big(\phi(a)p + \F\p(a)\big) =
  \phi(a)G(p) + \F\p(a) =
  \phi(a)h + \F\p(a).
  $$
  On the other hand,  given any $h$ satisfying \lcite{\haphihfp}, define
  $$
  G: ap+bq \in \X \mapsto  ah+b(1-h) \in A.
  $$
  Then, evidently $G$ is a left-$A$-linear map, which coincides with the identity on $A$.
Moreover, for any $x\in A$, we have
  $$
  %
  G\big((ap+bq)x\big) =
  G\big((a-b)px+bx\big) \$=  (a-b) G(px) + bx \={\paphiFp}
  (a-b)G\big(\phi(x)p+\F\p(x)\big)+bx \$=
  (a-b)\big(\phi(x)h+\F\p(x)\big)+bx \={\haphihfp}
  (a-b)hx+bx=
  G(ap+bq)x.
  $$
  Therefore $G$ is a conditional expectation.  Finally, it is clear that the
correspondence $G\mapsto h$ is injective.  \endProof

Condition \lcite{\haphihfp} is quite interesting because it relates
the intrinsic map $\F$ to the intrinsic
automorphism $\phi$.  Exploring this condition a bit further we obtain:

\state Proposition \label ExploreCondition
  Given any conditional expectation $G:\X\to A$, let $h=G(p)$ and $k = 1-h$.  Then,
  \izitem
  \halign{ \hfill \rm (\rzitemplus) \kern2pt  #& $#$ & = $#,$ \pilar{12pt}\cr
  & \E(a)   & ha +\phi(a)k \cr
  & \F(a)   & ka +\phi(a)h \cr
  & \E\p(a) & ka -\phi(a)k \cr
  & \F\p(a) & ha -\phi(a)h \cr}
  \medskip \noindent  for every $a\in A$.

\Proof
  Point (iv) follows immediately from \lcite{\haphihfp}.  As for (ii)
we have
  $$
  \F(a) = a-\F\p(a) \={iv}
  a - \big(ha-\phi(a)h\big) =
  (1-h)a +\phi(a)h =   ka +\phi(a)h,
  $$
  The proofs of (i) and (iii) follow similar lines.
\endProof

The following gives an important covariance condition relating $\phi$ and the images of $p$ and $q$ under a conditional expectation.

\state Proposition \label hkCovar Given a conditional expectation $G$ from $\X$ to $A$, let
  $h = G(p)$,
  and
  $k = G(q) = 1-h.$
  Then
  $$
  hka = \phi^2(a)hk \for a\in A.
  $$

\Proof
  We initially observe that
  $$
  \phi \E\p = -\F\p \E\p =  \F\p\phi,
  $$
  which may be easily proven by writing $\phi = \E-\F\p = \F-\E\p$.
  Secondly, by \lcite{\ExploreCondition.iii}  we have that
  $
  ka = \E\p(a) +\phi(a)k,
  $
  for every $a\in A$.
  Therefore
  $$
  hka =   h\big(\E\p(a) +\phi(a)k\big) \={\haphihfp}
  \phi\big(\E\p(a)\big)h + \F\p\big(\E\p(a)\big) +
  \phi^2(a)hk + \F\p\big(\phi(a)\big)k \$=
  -\F\p\big(\E\p(a)\big)h + \F\p\big(\E\p(a)\big) + \phi^2(a)hk - \F\p\big(\E\p(a)\big)k \$=
  -\F\p(\E\p(a))\big(h+k\big) + \F\p\big(\E\p(a)\big) + \phi^2(a)hk = \phi^2(a)hk.
  \endProof

In the presence of the automorphism $\Phi$ of \lcite{\EnterPsi}, it is interesting to
characterize the conditional expectations which commute with $\Phi$.

\state Proposition \label CondForCovar Let $G: \X\to A$ be a
conditional expectation and let $h=G(p)$.  Then the following are
equivalent
  \izitem
  \zitem $G\Phi = \Phi  G$,
  \zitem $\phi(h) = 1-h$.

\Proof
  Assuming (i) we have
  $$
  \phi(h) = \Phi(h) = \Phi(G(p)) = G(\Phi(p)) = G(q) = G(1-p) = 1-h.
  $$
  Conversely, if (ii) is known to hold, given any $c\in \X$, write $c=ap+bq$,  according to
\lcite{\UniqueDecomp} and notice that
  $$
  G\big(\Phi(c)\big) =
  G\big(\phi(b)p+\phi(a)q) = \phi(b)h+\phi(a)(1-h) \$=
  \phi(ah+b(1-h)) = \phi(G(c)) = \Phi(G(c)).
  \endProof

Given the relevance of conditional expectations satisfying the above
conditions we make the following:

\definition \label DefineCovar A conditional expectation $G:\X\to A$ satisfying the
equivalent conditions of \lcite{\CondForCovar} will be called \"{covariant}.

An interesting question is whether or not there exists a covariant conditional
expectation.  In the following sections we will give an affirmative answer in the context of
C*-algebras.

\section Strict C*-alloys with $B=\C^2$

From this point on we shall soup up the working hypothesis of the previous
two sections by assuming that we are given a unital strict C*-alloy
$\genBlend$ in which $B=\C^2$.

Since the present situation is a special case of the situation treated
above we may use all of the results so far obtained, including the
existence of the intrinsic maps.

Among the consequences of our strengthened standing hypothesis notice
that the projections $p$ and $q$ are now self-adjoint.

  As before we will refer to our alloy using the simplified notation  $(A,\C^2,\X)$.

As seen in \lcite{\FormulaForMult}, the left intrinsic pair encodes
the multiplication operation on $\X$.   In the present case we will
now show that the star operation may also be recovered from the left
intrinsic pair.

\state Proposition \label FormulaForStar Given $a, b\in A$, let
$c = ap+bq$.  Then
  $$
  c^* =   \big(\E(a^*) + \E\p(b^*)\big) p + \big(\F\p(a^*)+ \F(b^*)\big)q.
  $$

\Proof
  Left for the reader.
  \endProof

Let us now study the continuity of the intrinsic maps.

\state Proposition $\E$, $\F$ and  $\phi$ are bounded linear maps.

\Proof
  We first prove the continuity of $\E$ using the closed graph Theorem.  For this, assume
that $\{a_n\}_n$ is a sequence of elements in $A$, converging to some $a\in A$, and such that
$\E(a_n) \to b$.
  Then
  $$
  bp =
  \lim_{n\to\infty} \E(a_n)p =
  \lim_{n\to\infty} pa_np = pap = \E(a)p,
  $$
  so, $b=\E(a)$, and the continuity of $\E$ is established.   Similarly one proves that
$\F$ is continuous and hence so if $\phi$.
  \endProof

Other important properties relating the intrinsic maps and the metric
structure of the C*-algebras involved are discussed next.

\state Proposition \label SomeIneqs
  There exists a constant $K>0$, such that, for all $a$ in $A$, one
has that
  \izitem
  \zitem $\|ap\| \geq K\|a\|$,
  \zitem $\|aq\| \geq K\|a\|$,
  \zitem $\|\E(a^*a)\| \geq K^2\|a\|^2$,
  \zitem $\|\F(a^*a)\| \geq K^2\|a\|^2$,

\Proof  Observing that $\X = Ap+Aq$, one concludes that
  $$
  Ap = \X p = \{c\in \X: cq=0\}.
  $$
  This implies that $Ap$ is closed in $\X$, and hence that it is a
Banach space.  Since the map
  $$
  \lambda: a\in A \mapsto ap \in Ap,
  $$
  is a continuous bijection, the open mapping Theorem implies that it
is bicontinuous.  The constant $K = \|\lambda\inv\|\inv$ may then be
shown to satisfy (i).

  A similar constant can be chosen satisfying (ii) and so we may
rename $K$ to be the smallest of the two, and both (i) and (ii) will be
satisfied with the same constant $K$.

Next notice that
  $$
  \|\E(a^*a)\| \geq
  \|\E(a^*a)p\| =
  \|pa^*ap\| =
  \|(ap)^*ap\| =
  \|ap\|^2 \geq K^2\|a\|^2,
  $$
  proving (iii), while (iv) follows by a similar reasoning.
\endProof

When dealing with C*-algebras, conditional expectations
  are always required to be positive.  When a conditional expectation
$G$ also satisfies
  $$
  G(xx^*) = 0 \ \imply \ x=0,
  $$
  we say that $G$ is \"{faithful}.

  The following will be helpful later.

\state Lemma \label ahzero Given a faithful conditional expectation $G: \X \to A$, let
  $$
  h=G(p) \and k=G(q)=1-h.
  $$
  Then, for every $a\in A$, either
  $ah=0$, or $ak=0$, imply that $a=0$.

  \Proof Assuming that $ah=0$, notice that
  $$
  G\big((ap)(ap)^*\big) =  G(apa^*) =   aG(p)a^* = aha^* = 0,
  $$
  so $ap=0$, and hence \lcite{\UniqueDecomp} applies to give $a=0$.  The same conclusion may
similarly be obtained if we assume that $ak=0$. \endProof

Let us now prove some useful properties of covariant conditional
expectations.  Recall from \lcite{\DefineCovar} that a conditional
expectation is said to be covariant when $\phi\big(G(p)\big) = 1-G(p)$.

\state Corollary \label InvertibleGuys For any covariant
conditional expectation
  $G:\X\to A$, one has that
  \izitem
  \zitem $G(p)$ is invertible,
  \zitem $1-G(p)$ is invertible,
  \zitem $G$ satisfies the Pimsner--Popa finite index
condition \cite{\PP},
  \zitem $G$ is faithful.

\Proof
  Let $h = G(p)$.   Since $G$ is positive we have that $h$ is
positive as well, and so the spectrum of $h$, here denoted $\sigma(h)$,  is contained in the interval
$[0,+\infty)$.  Arguing by contradiction, suppose that
$0\in\sigma(h)$.  Consider the real valued functions
  $$
  f_n: [0,+\infty) \to \R
  $$
  defined by
  $$
  f_n(t) = \left\{\matrix{ \stake{10pt}
    1 -nt, &\hbox{if } t\leq 1/n, \cr
    0    , &\hbox{otherwise.}}\right.
  $$
  It is then easy to see that the element
  $
  a_n := f_n(h)
  $
  is positive and
  satisfies
  $$
  \|a_n\| = 1 \and \|a_nh\|\leq{1\over 4n}.
  $$
  We then have
  $$
  \E(a_n) \={\ExploreCondition.i}
  ha_n +\phi(a_n)(1-h) \={\CondForCovar.ii}
  ha_n +\phi(a_nh),
  $$
  whence
  $$
  \|\E(a_n)\| \leq \|ha_n\| +\|\phi\| \|a_nh\| \leq
  {1+\|\phi\|\over 4n}.
  $$
  However,  plugging $a = a_n\half$ in \lcite{\SomeIneqs.iii} leads to
  $$
  \|\E(a_n)\| \geq K^2 \|a_n\half\|^2 = K^2 \|a_n \| = K^2,
  $$
  bringing about a contradiction.  It follows that
$0\notin\sigma(h)$, proving (i).  Similarly one proves (ii).

In order to prove (iii), observe that $1-h = G(q)$, and hence $1-h$ is also
positive.  Being invertible, we have that
  $$
  1-h \geq \alpha >0,
  $$
  for some real number $\alpha$.  Since the same analysis applies to
$h$, we may also assume that
  $h\geq \alpha $.

  Given $c = ap+bq\in \X$, we then have
  $$
  G(cc^*) =   G\big((ap+bq)(pa^*+qb^*)\big) =
  G(apa^*+bqb^*) \$= aha^* + b(1-h)b^* \geq
  \alpha (aa^* + bb^*) \geq   \alpha (apa^* + bqb^*) = \alpha cc^*.
  $$
  This proves (iii), and (iv) is then an obvious consequence.
\endProof


Let us now discuss the consequences of left-right symmety imposed by the existence of
the star operation.

\definition If $f$ is any linear map between two *-algebras, we let $f^*$ be given by
  $$
  f^*(x) = f(x^*)^*.
  $$

Evidently $f^*$ is also linear.  In addition $f^*$ coincides with $f$ if and only if $f$ is *-preserving.

\state Proposition \label StarOtherSide The components of the left and right intrinsic
pairs and the intrinsic automorphisms satisfy:
  \izitem
  \zitem $\E^*=\Er$,
  \zitem $\F^*=\Fr$,
  \zitem $\phi^* = \phi\inv$,
  \zitem $\Phi^*=\Phi\inv$.

\Proof Given $a$ in $A$, observe that
  $$
  p\Er(a) = pap = (pa^*p)^* = \big(\E(a^*)p\big)^* = p\E(a^*)^* = p\E^*(a),
  $$
  so we deduce from \lcite{\UniqueDecomp} that $\Er(a) = \E^*(a)$.  The proof that $\Fr(a) =
\F^*(a)$ is similar.  With respect to (iii), recall from \lcite{\InverseOfPhi} that
  $$
  \phi\inv = \phir = \Er + \Fr - \id_A = \E^* + \F^* - \id_A = \phi^*.
  $$
  Addressing the last point,  observe that,
  given $a,b\in A$,
  $$
  \Phi^*(ap+bq)=
  \big(\Phi(pa^*+qb^*)\big)^* =
  \big(q\phi(a^*)+p\phi(b^*)\big)^* =
  \phi^*(a)q+  \phi^*(b)p,
  $$
  from where it is evident that $\Phi\Phi^*$ is the identity map, hence concluding the proof.
\endProof

Recall from \cite[4.1]{\Oka} that if $\rho $ is an automorphism of a C*-algebra, its \"{dual}
$\rho'$ is defined by
  $$
  \rho' = (\rho^*)\inv.
  $$
  Clearly $\rho$ is a *-automorphism if and only if $\rho'\rho=\id$.

According to \cite[6.1]{\Oka}, $\rho$ is said to be self-dual if $\rho'=\rho$.  If, in
addition, its spectrum\fn{Here $\rho$ is viewed simply as a bounded linear transformation and
hence one may speak of its spectum in the usual way.} consists of non-negative real numbers, one says that $\rho$ is positive.

For the convenience of the reader we state below a slight variant of the ``{polar decomposition
for isomorphisms}" from \cite{\Oka}, which we shall use in the sequell.

\state Theorem {\rm \cite[7.1]{\Oka}} Any automorphism  $\rho$  of a C*-algebra is written uniquely
as $\rho=\pi\gamma$, where $\pi$ is a *-automorphism and $\gamma$ is a positive automorphism.

An application of these ideas gives us some important information.

\state Proposition \label InvolutionCommute The automorphism $\Phi$
introduced in \lcite{\EnterPsi} is self-dual.  In addition, if
  $$
  \Phi=\Pi\,\Gamma
  $$
  is the polar decomposition of\/ $\Phi$,  then $\Pi$ is an involution commuting  with $\Gamma$.

\Proof
  Since $\Phi\inv=\Phi^*$ by \lcite{\StarOtherSide.iv}, we deduce that $\Phi'=\Phi$.
  We then have that
  $$
  \Phi = \Phi' = (\Pi\,\Gamma)' = \Gamma'\Pi' = \Gamma\,\Pi\inv = \Pi\inv(\Pi\,\Gamma\,\Pi\inv).
  $$
  It is evident that $\Pi\inv$ is a *-automorphism and that $\Pi\,\Gamma\,\Pi\inv$ is positive.  By
the uniqueness part of \cite[7.1]{\Oka} we conclude that $\Pi\inv = \Pi$, and
$\Pi\,\Gamma\,\Pi\inv=\Gamma$, concluding the proof.
  \endProof

By construction (see \lcite{\EnterPsi}) we have that $\Phi$ preserves $A$ and interchanges $p$
and $q$.  These properties  are reflected in the components of the polar decomposition as
follows:

\state Proposition \label ReflectedInPolar
  \izitem
  \zitem $\Gamma(A) = A = \Pi(A)$,
  \zitem $\Gamma(p) = p$,  and $\Gamma(q)=q$,
  \zitem $\Pi(p) = q$,  and $\Pi(q)=p$.

\Proof
  By the discussion following \cite[6.3]{\Oka}, one has that $\Phi^2=\Phi'\Phi$ is positive.
Moreover $\Gamma$ is the ``square root" of $\Phi^2$ in the sense of the analytical functional
calculus, where by square root we mean the principal branch of the complex square root function
defined on the open right half-plane.

Using Runge's Theorem we may find a sequence $\{f_n\}_n$ of polynomials converging uniformly to
the above mentioned square root function over some open neighborhood of the spectrum of
$\Phi^2$.  It follows that
  $$
  \Gamma = \sqrt{\Phi^2} = \lim_{n\to\infty} f_n(\Phi^2).
  \subeqmark Polys
  $$
  Since $\Phi^2(A)=A$, we then conclude that $\Gamma(A)\subseteq A$.

Observing that $\Gamma\inv = \sqrt{\Phi^{-2}}$, the same reasoning above gives
$\Gamma\inv(A)\subseteq A$, hence proving that $\Gamma(A)=A$.  To conclude the proof of (i) it is
now enough to observe that
  $$
  \Pi(A) = \Phi\big(\Gamma\inv(A)\big) =  \Phi(A) = A.
  $$

Since $\Phi$ interchanges $p$ and $q$, as already observed, one has that $\Phi^2$ fixes $p$ and
$q$.  Using \lcite{\Polys} we then have that
  $$
  \Gamma(p) = \lim_{n\to\infty} f_n(\Phi^2)(p) = \lim_{n\to\infty} f_n(1)\ p = p,
  $$
  and likewise $\Gamma(q)=q$.  Consequently
  $$
  \Pi(p) = \Phi\big(\Gamma\inv(p)\big) = \Phi(p) = q,
  $$
  and similarly $\Pi(q)=p$.
 \endProof

Since $\phi$ is the restriction of $\Phi$ to $A$,  we may easily
obtain the polar decomposition of the former, knowing that of the
latter:

\state Proposition  \label SmallpiAndgamma Let $\pi$ and $\gamma$ be the restrictions of\/ $\Pi$
and $\Gamma$ to $A$, respectively.  Then $\pi$ and $\gamma$ are
automorphisms of $A$ and
  $$
  \phi = \pi\gamma
  $$
  is the polar decomposition of $\phi$.  Moreover $\pi$ is an
involution commuting with $\gamma$.

\Proof
  By \lcite{\ReflectedInPolar.i} we have that $A$ is invariant under
$\Pi$ and $\Gamma$ and hence $\pi$ and $\gamma$ are indeed
automorphisms of $A$.  By \cite[6.3]{\Oka} we have that $\gamma$ is
positive and it is evident that $\pi$ is a *-automorphism.  Since it
is clear that   $\phi = \pi\gamma$, the uniqueness of the polar
decomposition \cite[7.1]{\Oka} warrants that to be the polar
decomposition of $\phi$.
  The last part of the statement follows from
\lcite{\InvolutionCommute}.
\endProof

We may now prove the existence of covariant conditional expectations.

\state Theorem \label ExistCovarExp Suppose that there exists a faithful conditional
expectation $\hat G: \X \to A$.  Then there exists another faithful conditional expectation $G$
such that, setting $h = G(p)$, one has
  \izitem
  \zitem $\pi(h) = 1-h$,
  \zitem $\phi(h) = 1-h$,
  \medskip\noindent and therefore $G$ is covariant.

\Proof
  Since $\Pi$ is an involutive *-automorphism of $\X$ preserving $A$, it is clear that
$\Pi\hat G\Pi$ is another conditional expectation onto $A$, and so is the map $G$ defined by
  $$
  G(c) = {\hat G(c) + \Pi \hat G\Pi(c) \over 2} \for c\in \X.
  $$
  It is easy to see that $ G$ is also faithful, since, for all $c\in \X$, one has that
  $$
  0\leq \hat G(c^*c) \leq 2 \, G(c^*c).
  $$
  Letting
  $$
  h := G(p) \={\ReflectedInPolar.iii} {\hat G(p) + \pi \big(\hat G(q)\big) \over 2},
  $$
  we claim that $\pi(h) = 1-h$.
  In order to prove it
  notice that
  $$
  \pi(h) =
  {\pi\big(\hat G(p)\big) + \hat G(q) \over 2} =
  {\pi\big(\hat G(1-q)\big) + \hat G(1-p) \over 2} \$=
  1- {\pi\big(\hat G(q)\big) + \hat G(p) \over 2} = 1-h,
  $$
  proving (i).
  We next claim that $\phi^2(h) = h$.  Indeed, if $k = 1-h$, then
\lcite{\hkCovar} applies and hence
  $$
  hkh = \phi^2(h)hk,
  $$
  which says that
  $$
  \big(h- \phi^2(h)\big)hk = 0,
  $$
  so the claim follows from \lcite{\ahzero}.

Evidently we also have that $\Phi^2(h) = h$ so,
  if $\{f_n\}_n$ is the sequence of polynomials employed in the proof of
\lcite{\ReflectedInPolar}, we obtain
  $$
  \Gamma(h) = \lim_{n\to\infty} f_n(\Phi^2)(h) = \lim_{n\to\infty} f_n(1)\ h = h.
  $$
  It follows that
  $$
  \phi(h) = \Phi(h) = \Pi\big(\Gamma(h)\big) = \Pi(h) = \pi(h) = 1-h.
  $$
  The last statement follows from \lcite{\CondForCovar}.
\endProof

The existence of covariant conditional expectations allows us
to prove that the automorphism $\gamma$ of \lcite{\SmallpiAndgamma}
is inner, as follows:

\state Proposition \label InnerAuto Let $G$ be a covariant conditional
expectation and let $h=G(p)$ and $k=1-h$.  Then
  $$
  \gamma(a) =  (hk)\half a (hk)\mhalf
  \for a\in A.
  $$

\Proof
Since both $h$ and $k$ are invertible by \lcite{\InvertibleGuys} we may use \lcite{\hkCovar} to write
  $$
  (hk)a(hk)\inv =  \phi^2(a) = \phi'\phi(a) = \gamma^2(a),
  $$
  where we have used that the intrinsic automorphism $\phi$ is
self-dual.  The conclusion then follows from \cite[6.6]{\Oka}.
  \endProof

We now summarize our main results so far.

\state Theorem \label MainConcrete Let $(A,\C^2,\X)$ be a unital strict
C*-alloy such that there exists a faithful conditional expectation
$\hat G:C\to A$.  Then there exists a *-automorphism $\pi$ of $A$, and
a positive invertible element $h$ in $A$, such that
  \izitem
  \zitem $\pi^2=\id$,
  \zitem $\pi(h) = 1-h$.
  \medskip\noindent
  Moreover,  letting $k = 1-h$, one has that the intrinsic
automorphism and the components of the left intrinsic pair are entirely
determined in terms of $\pi$ and $h$ by
  $$ \def\a{\hfill\pilar{14pt}\cr}
  \matrix{
  \phi(a) &=& \pi\big((hk)\half a (hk)\mhalf\big) &=& (hk)\half \ \pi(a) \ (hk)\mhalf, \a
  \E(a) &=& ha + \phi(a)k, \a
  \F(a) &=& ka + \phi(a)h, \a
  }
  $$
  for all $a$ in $A$.
  Finally, the map
  $$
  G: ap+bq \in \X \mapsto ah+bk\in A
  \subeqmark CanonicCondExp
  $$
  is a faithful covariant conditional expectation satisfying
the Pimsner--Popa finite index condition.

\Proof
  Let $G$ be the faithful covariant conditional expectation given by
\lcite{\ExistCovarExp}, let $h = G(p)$, and let
$k=1-h=G(q)$. Therefore
  $$
  G(ap+bq) = ah+bk \for a,b\in A,
  $$
  and the last statement is then a consequence of \lcite{\InvertibleGuys}.

  Since $p$ is positive and $G$ preserves positivity, it is clear that
$h$ is positive.  By \lcite{\InvertibleGuys} we have that $h$ and $k$ are invertible.

Letting $\pi$ be as in \lcite{\SmallpiAndgamma},  one sees that (i) follows
from the last statement in \lcite{\SmallpiAndgamma},  while (ii) is the content of \lcite{\ExistCovarExp.i}.

The first formula above for $\phi$ is a consequence of
\lcite{\SmallpiAndgamma} and \lcite{\InnerAuto}, while the second one
follows from the fact that $hk$ is a fixed point for $\pi$, as one may
easily verify using (ii).  The expressions for $\E$ and $\F$ above in
turn follow from \lcite{\ExploreCondition.i--ii}.
  \endProof

The relevance of the above result is that the whole algebraic and analytical
structure of $\X$ may be recovered from $\pi$ and $h$.  This is because, once
the intrinsic pair is known, we may use \lcite{\FormulaForMult} and
\lcite{\FormulaForStar} to recover both the multiplication and star
operations on $\X$.  The norm may also be recovered
since the norm on any C*-algebra is encoded in its
*-algebraic structure:  the norm of an element $x$ coincides with the
square root of the spectral radius of $x^*x$.  It is
therefore meaningful to give the folowing:

\definition We will say that the pair $(\pi,h)$ is the \"{fundamental data}
for the alloy $(A,\C^2,\X)$.

\section Strict C*-alloys and crossed products by $\Z_2$

As seen in the previos section, the whole structure of a unital strict C*-alloy may
be recovered from its fundamental data.  A natural question arising
from this is whether or not one may construct a C*-alloy from a pair
$(\pi, h)$, where $\pi$ is a *-automorphism of a unital C*-algebra $A$, and
$h$ is a positive invertible element in $A$ satisfying
\lcite{\MainConcrete.i-ii}.

This may in fact be done by first
\"{defining} operators $\E$ and $\F$ on $A$ using the formulas
provided by \lcite{\MainConcrete}.  On the vector space $\X = A\oplus
A$, where a pair $(a, b)$ is formally denoted $ap+bq$, we may then
introduce a *-algebra structure by employing formulas
\lcite{\FormulaForMult} and \lcite{\FormulaForStar}.  It may then be
shown that ${(A,\C^2,i,j,\X)}$ is a strict C*-alloy, where
  $$
  \matrix{
  i & : & a\in A & \mapsto & (a,a)\in \X, & \hbox{and}\cr \pilar{15pt}
  j & : & (\lambda, \mu)\in \C^2 & \mapsto & (\lambda,\mu)\in \X.}
  $$

The situation is however even more interesting in the sense that $\pi$
alone provides enough information to construct $\X$,  while $h$ is necessary only in
order to locate
$p$ and $q$ within $\X$.  In order to explain this in detail let us
first quickly analyze the C*-alloy  arising from an action of $\Z_2$.

If $A$ is a unital C*-algebra and $\pi$ is an involutive
*-automorphism of $A$ we may form an action of $\Z_2$ on $A$ by
mapping
the generator of $\Z_2$ to $\pi$.
  We may therefore consider the crossed product algebra
  $$
  \X = A\rtimes_\pi \Z_2.
  $$
  The unitary element implementing the action of $\pi$, here denoted
by $\varpi$, is clearly a self-adjoint unitary.   Therefore,
  defining
  $$
  p = {1+\varpi \over 2}
  \and
  q = {1-\varpi \over 2},
  $$
  we have that $p$ and $q$ are complementary projections which
therefore generate a copy of $\C^2$ within $\X$.  Thinking of $A$ as a
subalgebra of $\X$ in the usual way, observe that, for every $a$ in
$A$, one has that
  $$
  pap = {1\over 4}(1+\varpi )a(1+\varpi ) =
  {1\over 4}(a+\varpi a+a\varpi +\varpi a\varpi ) \$=
  {1\over 4}\big(a+\pi(a)\varpi +a\varpi +\pi(a)\big) =
  {a+\pi(a)\over 2} \, {1+\varpi \over 2} =
  \E(a)p,
  $$
  where $\E(a)$ is defined to be $(a+\pi(a))/2$.  Routine
calculations show, in fact, that $(A,\C^2,\X)$ is a strict C*-alloy
with left intrinsic pair $(\E, \E)$, meaning that $\F=\E$.
  The intrinsic automorphism,  namely $\phi=\E+\F-\id$,  therefore
coincides with $\pi$.

Recall that the usual conditional expectation $G:\X\to A$ is given by
  $$
  G(a+b\varpi ) = a
  \for a,b\in A.
  $$
  The elements $h$ and $k$, which played crucial roles above, are therefore
given by
  $$
  h = G(p) = {1\over 2}
  \and k = 1-h = {1\over 2}.
  $$

We have thus proven:

\state Proposition If $\pi$ is an involutive *-automorphism of a unital
C*-algebra $A$, and $\X=A\rtimes_\pi\Z_2$, then  $(A,\C^2,\X)$ is a
strict C*-alloy.  Moreover its intrinsic automorphism coincides with $\pi$,
and $h=1/2$.

In our next result we shall prove that, under suitable
conditions, the above example is essentially the only one.

\state Theorem
  Let $(A, \C^2, \X)$ be a unital strict C*-alloy.  Suppose moreover that there exists a
faithful conditional expectation $\hat G:\X\to A$, and let $(\pi, h)$ be
the fundamental data provided by \lcite{\MainConcrete}.
  Then there exists
an involutive *-automorphism $\pi$ of $A$ and a *-isomorphism
  $$
  \rho : A\rtimes_\pi\Z_2 \to \X,
  $$
  coinciding with the identity on $A$, and such that
  $$
  \matrix{
  \rho(\hfill\varpi\hfill ) &=& (hk)\mhalf(kp-hq), & \hbox{and}\cr \pilar{15pt}
  \rho\big((hk)\half \varpi  +  h\big) &=& p,
  }
  $$
  where $\varpi $ is the unitary element implementing the action of
$\pi$, and $k=1-h$.

  \newcount \claimno \claimno = 0
  \def\Claim#1\par{\global \advance \claimno by 1 \bigskip\noindent {\tensc Claim
\number \claimno: }#1\bigskip}

\Proof
  By direct computation,  using the expression for $\phi$ given in
\lcite{\MainConcrete}, one may show that
  $$
  \phi(h)=k \and \phi(k)=h.
  \subeqmark PhiExchangeshk
  $$

From mow on   the proof will consist of a series of claims, starting with:

\Claim $\E(hk) = hk = \F(hk)$,  and $\E\p(hk) = 0 = \F\p(hk)$.

In order to verify this we compute
  $$
  \E(hk) \={\ExploreCondition.i} h(hk) + \phi(hk)k \={\PhiExchangeshk}
  hhk + khk = (h+k)hk = hk,
  $$
  while
  $$
  \F(hk) \={\ExploreCondition.ii} k(hk) + \phi(hk)h \={\PhiExchangeshk}
  khk + khh = kh(k+h) = hk.
  $$

\Claim  $\E(h\inv k) = 1 = \F(k\inv h)$.

We have
  $$
  \E(h\inv k) \={\ExploreCondition.i}
  h(h\inv k) + \phi(h\inv k)k \={\PhiExchangeshk}
  h(h\inv k) + (k\inv h)k = k+h = 1,
  $$
  while
  $$
  \F(k\inv h) \={\ExploreCondition.ii}
  k(k\inv h) + \phi(k\inv h)h \={\PhiExchangeshk}
  k(k\inv h) + (h\inv k)h = h+k = 1.
  $$

We will now introduce the element of $\X$ which will correspond to the
implementing unitary in the crossed product.  Let
  $$
  u = (hk)\mhalf(kp-hq).
  $$
  It is our next immediate goal to show that $u$ is a self-adjoint
unitary.  We begin with the following:

\Claim $u$ is an isometry.

  We have
  $$
  u^*u =
  (kp-hq)^*(hk)\inv(kp-hq) =
  (pk-qh)(hk)\inv(kp-hq) \$=
  (ph\inv-qk\inv)(kp-hq) =
  ph\inv kp+qk\inv hq \$=
  \E(h\inv k)p+\F(k\inv h)q =
  p+q=1.
  $$

\Claim $hk$ commutes with $p$ and $q$.

  Observe that
  $$
  phk \={\EnterEnadF.i}
  \E(hk)p + \F\p(hk)q =
  hkp.
  $$
  Since $q=1-p$,  it is clear that $hk$ also commutes with $q$.

Our next claim refers to the second factor of $u$, namely the element
defined by
  $$
  v=kp-hq.
  $$

\Claim One has that $v = p-h$,  and hence $v$ is self-adjoint.

  Notice that
  $$
  v = kp - hq = (1-h)p - hq = p - hp - hq = p-h(p+q) = p-h.
  $$
  Since both $p$ and $h$ are self-adjoint, the claim follows.

\Claim $u$ is self-adjoint.

Since $hk$ commutes with $p$ and $q$, it is clear that $hk$ also
commutes with $v$.  Therefore $(hk)\mhalf$ commutes with $v$.  The
claim is then a consequence of the obvious fact that the product of
two commuting self-adjoint elements is again self-adjoint.

\Claim $u$ is unitary.

We have already seen that $u^*u=1$,  so the result follows from the fact
that $u=u^*$.

\Claim For every $a$ in $A$ one has that $\pi(a) = uau\inv$.

Recalling from claim \lcite{5} that $v = p-h$, one has
  $$
  va = pa-ha \={\haphihfp\ \& \paphiFp}
  \phi(a)p+\F\p(a) - \phi(a)h-\F\p(a) =
  \phi(a)(p-h) =   \phi(a)v.
  $$
  Since $v = (hk)\half u$, one sees that $v$ is invertible, and so
  $$
  \phi(a) = vav\inv.
  $$
  Consequently
  $$
  \pi(a) = (hk)\mhalf\phi(a)(hk)\half =
  (hk)\mhalf vav\inv(hk)\half = uau\inv.
  $$

By the universal property of crossed products there is a
*-homomorphism
  $$
  \rho : A\rtimes_\pi\Z_2 \to \X,
  $$
  extending the identity map on $A$, and sending the implementing
unitary $\varpi $ to $u$.

We then have that
  $$
  \rho\big((hk)\half \varpi + h\big) =
  (hk)\half \rho(\varpi) + h =
  (hk)\half u + h =
  v + h = p.
  $$
  It now remains to show that $\rho$ is bijective.
  By the computation above we see that $p$ is in the range of $\rho$,
and this in turn shows that $\rho$ is onto.

  Let $G$ be the faithful conditional expectation provided by
\lcite{\MainConcrete}.  Then $G(p)=h$, and hence
  $$
  G(u) = G\big((hk)\mhalf(p-h)\big) = (hk)\mhalf \big(G(p)-h\big) = 0.
  $$

  If $H$ is the standard conditional expectation from
$A\rtimes_\pi\Z_2$ to $A$, namely that which is given by
  $$
  H(a+b\varpi) = a
  \for a,b\in A,
  $$
  we claim that $G\rho = \rho H$.  In fact, we have
  $$
  G\big(\rho(a+b\varpi)\big) =
  G(a+b u) = a+b G(u) = a = \rho\big(H(a+b\varpi)\big).
  $$

In order to prove that $\rho$ is one-to-one,  assume that $x\in
A\rtimes_\pi\Z_2$ is such that $\rho(x)=0$.  Then
  $$
  0 = G\big(\rho(xx^*)\big) =
  \rho\big(H(xx^*)\big).
  $$
  Since $H(xx^*)$ lies in $A$, and since $\rho$ coincides with the
identity on $A$, and hence is injective there,  we deduce that
$H(xx^*)=0$.  Finally, since $H$ is faithful, we conclude that $x=0$,
thus proving that $\rho$ is injective.
\endProof

\section A non strict finite dimensional example

\label BadExample When dealing with tensor products, one is used to
believe that most analytical problems disappear provided one of the
factors is finite dimensional.

Contrary to such expectations, in this section we shall give an
example a unital C*-blend $(A, B, \X)$, with $B$ being a finite dimensional
algebra, which is not strict.  In particular, upon identifying $A\*B$
with a subspace of $\X$, we may view $A\*B$ as a normed space which
will turn out not to be complete,  regardless of the fact that $B$ is
finite dimensional.

In order to prepare for the construction of our counter-example we
first consider the following elementary example:
  let $\X$ denote the C*-algebra formed by all $2\times2$ complex matrices.  Given
any real number
  $r$ in the open interval $(0,1)$, consider the element of $\X$ given by
  $$
  p(r) =  \m
  {r }
  {\sqrt{r-r^2}\, }
  {\sqrt{r-r^2}}
  {1-r\pilar{15pt} },
  $$
  which is easily seen to be a projection.  We shall also consider its
complementary projection
  $$
  q(r)=1-p(r).
  $$
  We shall now describe two subalgebras of $\X$ which will form a blend
of C*-algebras.  On the one hand we will let $A$ be the subalgebra
formed by all diagonal matrices and, on the other, $B$ will be the
subalgebra generated by $p(r)$ and $q(r)$.  Clearly $B$ is isomorphic
to $\C^2$, so we will identify $B$ and $\C^2$ from now on (however we
will not give much attention to the fact that $A$ is also isomorphic
to $\C^2$).

  We leave it for the reader to prove that the triple $(A,B,\X)$,
henceforth also
referred to as $(A,\C^2,\X)$, is a C*-blend. It is also easy to see
that, if $a = \d xy\in A$, then
  $$
  \E(a) = \big(rx+(1-r)y\big) I_2
  \and
  \F(a) = \big((1-r)x+ry\big) I_2,
  \eqmark EnFnPn
  $$
  where $I_2$ is the identity $2\times2$ matrix.
  Consider the C*-algebra
  $$
  \X\f := \ell_\infty(\X)
  $$
  formed by all bounded sequences of elements in $\X$,  under pointwise
operations.   It is evident that
  $$
  A\f := \ell_\infty(A)
  $$
  sits as a subalgebra of $\X\f$.

\state Proposition Given any sequence $\{r_n\}_n$ of real numbers in the open interval
$(0,1)$, let $p$ and $q$ be the elements of $\X\f$ given by
  $$
  P = \big(p(r_n)\big)_{n\in\N} \and Q=1-P.
  $$
  Then $A\f P + A\f Q$ is a *-subalgebra of $\X\f$,  containing $A\f$.

\Proof
  Throuout this proof we will write $p_n$ for $p(r_n)$. Moreover, when $r$ is
replaced by $r_n$, the maps described in \lcite{\EnFnPn} will be
written $\E_n$ and $\F_n$, respectively.

  We first claim that
  $$
  PA\f\subseteq A\f P + A\f Q.
  $$
  In order to see this, let $a \in A\f$.  Then
  $$
  Pa = \vec{p_na_n} =
  \vec{\E_n(a_n)p_n + \F_n\p(a_n)q_n} =
  \vec{\E_n(a_n)}P + \vec{\F_n\p(a_n)}Q.
  $$
  Noticing that the $\E_n$ and the  $\F_n$  are uniformly bounded, the
above calculation implies the claim.

It is clear that $A\f$ is contained in  $A\f P + A\f Q$, and hence also
  $$
  QA\f = (1-P)A\f \subseteq A\f + PA\f \subseteq A\f P + A\f Q.
  $$
  It follows that
  $$
  PA\f+QA\f  \subseteq A\f P + A\f Q.
  $$
  Taking the adjoint on both sides above, we deduce  the
reverse inclusion and hence that
  $$
  PA\f+QA\f = A\f P + A\f Q.
  $$
  In particular it follows that $A\f P + A\f Q$
  is a self-adjoint set and it is now easy to prove all remaining
assertions.  \endProof

Based on the above result, the closure of $A\f P + A\f Q$,  here
denoted by
  $$
  \tilde \X = \overline{A\f P + A\f Q},
  $$
  is seen to be a
C*-algebra and   consequently
  $$
  (A\f,\C^2,\tilde \X)
  \eqmark BadBlend
  $$
  is a C*-blend, where we identify $\C^2$ with the
subalgebra of $\tilde \X$ spanned by $P$ and $Q$.

The main question we wish to address here is related to whether or not $A\f P + A\f Q$ is
closed.  In fact we wish to prove that this is not the case when the
$r_n$ tend to zero.

\state Proposition Let $\{r_n\}_n$ be a sequence in $(0,1)$ such that
  $$
  \lim_{n\to\infty}r_n=0.
  $$
  Then $A\f P + A\f Q$ is not closed in $\tilde \X$ and hence
$(A\f,\C^2,\tilde \X)$ is a C*-blend which is not strict.

\Proof
  We shall suppose, by way of contradiction, that $A\f P + A\f Q$
is closed,  and hence that $\tilde \X = A\f P + A\f Q$.

From this it follows that
  $$
  A\f P = \tilde \X P = \{x\in \tilde \X : xQ=0\},
  $$
  which is therefore closed in $\tilde \X$.  The map
  $$
  a\in A\f \to aP \in A\f P
  $$
  is clearly continuous and bijective, and hence, by the Open Mapping
Theorem, we deduce that it is bounded by below, that is, there exists
a constant $K>0$, such that
  $$
  \|aP\|\geq K\|a\|
  \for a\in A\f.
  $$
  Given a positive integer $m$, let $a=\vec{a_n}$ be the element of $A\f$ given by
  $$
  a_n = \delta_{n, m}\d 10.
  $$
  Then $aP$ has a single nonzero coordinate in the $m^{th}$ position,
and that coordinate is given by
  $$
  a_mp_m =
  \d 10 \m {r_m} {\sqrt{r_m-r_m^2}\, } {\sqrt{r_m-r_m^2}} {1-r_m\pilar{15pt} } =
  \m {r_m} {\sqrt{r_m-r_m^2}\, } {0}{0}.
  $$
  It follows that
  $$
  0<K = K\|a\| \leq \|aP\| = \|a_mp_m\| =
  \sqrt{r_m} \ \ {\buildrel m\to\infty \over \longrightarrow}\ 0,
  $$
  a contradiction. \endProof

\section C*-blends and commuting squares

Let us now explore a different class of C*-blends, not necessarily
coming from group actions.  The initial data we shall use in the
construction to be described in this section is a commuting square \cite{\Po}, by
which we mean that we are given C*-algebras $A$, $B$, $C$ and $D$,
such that
  $$
  \matrix{
  A & \supset & B \cr
  \vrule height 14pt depth 8pt width 0pt  \cup && \cup \cr
  C & \supset & D,
  }
  $$
  and $D=B\cap C$.  One is moreover given conditional expectations
  $$
  E:A\to B \and F:A\to C,
  $$
  satisfying $EF=FE$.   Evidently one then has that
  $
  G := E F
  $
  is a conditional expectation onto $D$.


Let us now quickly describe
  Watatani's version \cite{\W} of
  the celebrated Jones' basic construction relative to $G$.  One first
introduces a $D$-valued inner product on $A$ by the formula
  $$
  \langle a,b\rangle = G(a^*b), \quad \forall a,b\in A.
  \eqmark DefineInnProd
  $$
  The completion of $A$ relative to the norm
arising from this inner product is a  right Hilbert $D$-module, which
we shall denote by $M$.  The left action of $A$ on itself may be shown
to extend to a *-homomorphism
  $$
  \lambda : A \to \L(M),
  $$
  where $\L(M)$ denotes the C*-algebra of all adjointable operators on
$M$.
  The fact that
  $$
  G(a)^*G(a) \leq G(a^*a) \for a\in A,
  $$
  which may be easily proven by noting that $G(x^*x)\geq 0$, where $x =
a-G(a)$, implies that $G$ extends to a bounded linear operator on $M$,
which we denote by $g$,  and which may be shown to be a projection in
$\L(M)$,  often referred to as the \"{Jones projection}.
It is easy to show that
  $$
  g\lambda(a)g = \lambda\big(G(a)\big)g \for a\in A.
  \eqmark JonesRel
  $$

One may similarly show that $E$ and $F$
extend to bounded operators on $M$, respectively denoted $e$ and $f$,
providing projections $e$ and $f$ in $\L(M)$ satisfying
  $$
  e\lambda(a)e = \lambda(E(a))e
  \and
  f\lambda(a)f = \lambda(F(a))f
  \for a\in A.
  \eqmark MoreJonesRel
  $$
  The fact that $G=EF=FE$ immediately implies that
  $$
  g = ef=fe.
  \eqmark gEQUALef
  $$

  Given $a_1,b_1,a_2, b_2\in A$, observe that
  $$
  \big(\lambda(a_1)g\lambda(b_1)\big)\big(\lambda(a_2)g\lambda(b_2)\big) = \lambda\big(a_1G(b_1a_2)\big)g\lambda(b_2),
  $$
  which implies that the linear span of $\lambda(A)g\lambda(A)$ is an algebra, easily
seen to be self adjoint.  We thus obtain a closed *-subalgebra of
$\L(M)$ by setting
  $$
  K_g = \clspan{\lambda(A)g\lambda(A)}.
  $$
  In an entirely similar fashion we have the closed *-subalgebras
  $$
  K_e  = \clspan{\lambda(A)e\lambda(A)} \and   K_f  = \clspan{\lambda(A)f\lambda(A)}.
  $$

  \state Proposition
  One has that $K_eK_g$, $K_gK_e$, $K_fK_g$ and $K_gK_f$ are all
contained in $K_g$, and hence both $K_e$ and $K_f$ may be naturally
mapped into the multiplier algebra $M(K_g)$. We denote these maps
by
  $$
  i_e : K_e \to M(K_g) \and   i_f : K_f \to M(K_g).
  $$

\Proof  Given $a,b,c,d\in A$, and using that $g = eg$, we have that
  $$
  \big(\lambda(a)e\lambda(b)\big)\big(\lambda(c)g\lambda(d)\big) =
  \lambda(a)e\lambda(bc)eg\lambda(d) =
  \lambda\big(aE(bc)\big)g\lambda(d) \in K_g,
  $$
  which proves that $K_eK_g\subseteq K_g$.  The remaining inclusions
may be proven similarly.
  \endProof

\state Question \label MyQuestion \rm Is $\big(K_e,K_f,i_e,i_f,K_g\big)$ a C*-blend?

My interest in the whole idea of C*-blends actually arouse from this
question,  to which I still do not have a definitive answer in its full
generality.  However we at least have:

\state Proposition  \label OneContainement If either $B$ or $C$ contain an approximate
identity for $A$, then the ranges of $i\vezes j$ and $j\vezes i$ are
both dense in $K_g$.

\Proof
  By taking adjoints, it is enough to prove that
  $$
  K_g\subseteq \clspan{K_eK_f}.
  $$
  Let $\{u_i\}_i$ be an approximate
identity for $A$ contained, say, in $B$.  Then, given $x$ in the dense
image of $A$ within $M$, we
  have that
  $$
  e\lambda(u_i) f(x) = E(u_iF(x)) = u_iE(F(x)) = u_i G(x) = \lambda(u_i)g(x),
  $$
  which says that
  $e\lambda(u_i)f = \lambda(u_i)g$.  Therefore, given $a,b\in A$, we
have
  $$
  \lambda(a)g\lambda(b) =
  \lim_i\lambda(a)\lambda(u_i)g\lambda(b) =
  \lim_i\lambda(a) e\lambda(u_i)f\lambda(b) \$=
  \lim_i\big(\lambda(a) e\lambda(u_i\half)\big)
  \big(\lambda(u_i\half) f\lambda(b)\big)  \in \clspan{K_eK_f}.
\endProof

This may now be used to state a sufficient condition for \lcite{\MyQuestion}.

\state Proposition \label SufCondForBlend Suppose that either $B$ or
$C$ contain an approximate identity for $A$.  If
  $$
  e\lambda(A)f\subseteq K_g,
  $$
  then $\big(K_e,K_f,i_e,i_f,K_g\big)$ is a C*-blend.

\Proof
  Once in possession \lcite{\OneContainement} it is enough to prove
that the range of $i_e\vezes i_f$, and consequently also of
$i_f\vezes i_e$, is contained in $K_g$.  In other words we must show
that
  $$
  K_eK_f\subseteq K_g.
  $$
  Given $a,b,c,d\in A$, we have
  $$
  \big(\lambda(a)e \lambda(b)\big) \big(\lambda(c)f \lambda(d)\big) =
  \lambda(a)e \lambda(bc)f \lambda(d) \subseteq
  \lambda(a)  K_g \lambda(d) \subseteq K_g.
  \endProof

This is as much as we can say in the present generality, so let us now
consider a somewhat restrictive special case.
Recall from \cite{\W}, that the conditional expectation $G$ is said to
be of \"{index finite type} provided there is a finite set
  $$
  \{u_1,\ldots,u_n\} \subseteq A,
  $$
  called a \"{quasi-basis}, such that
  $$
  a = \sum_{i=1}^n u_iG(u_i^*a)
  \for a\in A.
  \eqmark FiniteIndCond
  $$

  Let us assume for the time being that $G$ is of index finite type.
As a consequence $G$ is necessarily faithful \cite[2.1.5]{\W} and,  as
one may easily show, $\lambda$ is injective.  We shall therefore
identify $A$ with its image under $\lambda$ without further warnings.

Some other aspects of the present situation are also much simplified because
of the following:

\state Proposition If $G$ is of index finite
type, let $\{u_1,\ldots,u_n\}$ be a quasi-basis.  Then
  \izitem
  \zitem $\sum_{i=1}^n u_igu_i^* = 1$,
  \zitem $\span{AgA} = K_g = \L(M) = L_D(M)$,  where the latter refers to the
set of all $D$-linear maps on $M$.

\Proof
  Initially notice that,  by \cite[2.1.5]{\W},
$A$ is already complete with the norm arising
from \lcite{\DefineInnProd}, so $M=A$.
  For every $a\in M$, one then has that
  $$
  \sum_{i=1}^n u_igu_i^*(a)=
  \sum_{i=1}^n u_iG(u_i^*a)= a,
  $$
  proving (i).
  As for (ii), we first notice that the inclusions
  $$
  \span{AgA} \subseteq K_g \subseteq \L(M) \subseteq L_D(M)
  $$
  are all evident.
  Given a $D$-linear map $T$ on $M$, notice that, for every $a\in M$,
  $$
  T(a) =
  T\Big(\sum_{i=1}^n u_iG(u_i^*a)\Big) =
  \sum_{i=1}^n T(u_i)G(u_i^*a)=
  \sum_{i=1}^n T(u_i)gu_i^*a,
  $$
  so we  see that
  $$
  T =
  \sum_{i=1}^n T(u_i)gu_i^* \in \span{AgA}.
  \endProof

Evidently both $e$ and $f$ are $D$-linear, as are all operators on
$M$ in sight, so we have that $eAf \subseteq K_g$.   Employing
\lcite{\SufCondForBlend} we then have the following affirmative answer
to \lcite{\MyQuestion}:

  \state Proposition If $G$ is of index finite type then $\big(K_e,K_f,i_e,i_f,K_g\big)$ is a C*-blend.

Let us now give another partial answer to \lcite{\MyQuestion}, this
time without assuming finiteness of the index.
  Still under the assumption that we are given a commuting square as
above, we will now assume the following:

\state {Standing Hypothesis} \label StandHypCommut $A$ is a unital commutative algebra and $D = \C\, 1$.

As a consequence $G$ is necessarily of the form
  $$
  G(a)=\phi(a)1 \for a\in A,
  \eqmark GisPhi
  $$
  where $\phi$ is a state of $A$.

Being a right Hilbert $D$-module, $M$ is nothing but a Hilbert space
while $\lambda$ is just the GNS representation associated to $\phi$.
Using $\lambda$ we will view $M$ as a left $A$ module, and hence we
will adopt the notation
  $$
  a\eta := \lambda(a)\eta
  \for a\in A, \ \eta\in M.
  $$
  Denoting the image of $1$ in $M$
by $\xi$, we evidently have the well known GNS formula
  $$
  \phi(a) = \<a\xi,\xi\>
  \for a\in A.
  $$

Let us now prove a crucial inequality to be used later.

\state Lemma \label MainIneq
  Given $a\in A$, and $c_1,\ldots,c_n\in C$, one has that
  $$
  \soma \|eac_i\xi\|^2 \leq \|E(a^*a)\|\ \|\mu\|,
  $$
  where $\mu$ is the $n\times n$ scalar matrix with $\mu_{ij}
= \<c_i\xi,c_j\xi\>$, for $i,j=1, \ldots, n$.

\Proof
  Since $E(1)=1$,  one has that $e(\xi)=\xi$,  and hence for every
$a\in A$ we have
  $$
  ea(\xi) = eae(\xi) = E(a)e(\xi) = E(a)\xi,
  $$
  and similarly
  $$
  fa(\xi) = F(a)\xi \for a\in A.
  \subeqmark RangeOfF
  $$
  We then have
  $$
  \soma \|eac_i\xi\|^2 =
  \soma \|E(ac_i)\xi\|^2 =
  \soma \phi\big(E(c_i^*a^*)E(ac_i)\big) =
  \phi(b),
  \subeqmark Somaeac
  $$
  where
  $
  b :=\soma E(c_i^*a^*)E(ac_i).
  $
  \ Let $m$ be the $n\times n$ matrix over $A$ given by
  $$
  m = \pmatrix{
  c_1 & 0 & \ldots & 0 \cr
  c_2 & 0 & \ldots & 0 \cr
  \vdots & \vdots & \ddots & \vdots\cr
  c_n & 0 & \ldots & 0 \cr
  }.
  $$
  Observing that
  $$
  E(m^*a^*)E(am) =
  \pmatrix{
  b & 0 & \ldots & 0 \cr
  0 & 0 & \ldots & 0 \cr
  \vdots & \vdots & \ddots & \vdots\cr
  0 & 0 & \ldots & 0 \cr
  },
  $$
  we see that
  $$
  \|b\| =
  \|E(m^*a^*)E(am)\| =
  \|E(am)\|^2.
  $$

  \bigskip Recall from \cite[2.9]{\Rieffel} that
  $$
  \def\x{x}\def\y{y}
  \|E(\x^*\y)\|^2 \leq \|E(\x^*\x)\|\ \|E(\y^*\y)\|
  \for \x,\y\in A,
  $$
  so we have
  $$
  \def\x{m^*} \def\xs{m} \def\y{a} \def\ys{\y^*}
  \|b\| =  
  \|E(am)\|^2 \={!}
  \|E(ma)\|^2 \leq
  \|E(\xs\x)\|\, \|E(\ys\y)\|.
  \subeqmark MMStar
  $$

  In case the reader is interested in attempting to eliminate the
commutativity hypothesis assumed in \lcite{\StandHypCommut}, we have marked the only two
places in which we used the commutativity of $A$ with an exclamation
mark, the first one appearing in the calculation just above.

  The entry  $(i,j)$ of the matrix $E(mm^*)$ is given by
  $$
  E(c_ic^*_j) =
  EF(c_ic^*_j) =
  \phi(c_ic^*_j) \={!}
  \phi(c^*_jc_i) =
  \<c_j^*c_i\xi,\xi\> =
  \<c_i\xi,c_j\xi\>,
  $$
  which precisely says that $E(mm^*)$ coincides with the matrix
$\mu$ in the statement.  We therefore have
  $$
  \soma \|eac_i\xi\|^2 \ \={\Somaeac}\
  \phi(b) \leq
  \|b\| \buildrel {(\MMStar)}
\over \leq
  \|\mu\|\|E(a^*a)\|.
  \endProof

Recall that $\xi$ is the image of 1 in $M$, and hence that
$g(\xi)=\xi$.  Moreover, given $a,b, c\in A$, let
$\eta=c\xi$, and notice that
  $$
  (agb^*)(\eta) =
  (agb^*)(c\xi) = agb^*cg\xi = aG(b^*c)g\xi = aG(b^*c)\xi \$=
  \phi(b^*c)a\xi = \<c\xi,b\xi\>a\xi =
  \<\eta,b\xi\>a\xi.
  $$
  Since the set of vectors of the form $\eta = c\xi$ is dense in $A$, this
shows that
  $$
  (agb^*)(\eta) =  \<\eta,b\xi\>a\xi,
  \eqmark MyRankOne
  $$
  for every $\eta$ in $M$.
  Given any two vectors $\zeta_1$ and $\zeta_2$ in
$M$, consider the rank-one operator
  $\Omega_{\zeta_1,\zeta_2}$, defined by
  $$
  \Omega_{\zeta_1,\zeta_2}(\eta) =  \<\eta,\zeta_2\>\zeta_1
  \for \eta\in M.
  $$
  By \lcite{\MyRankOne}, we then have that
  $$
  \lambda(a)g\lambda(b^*) =   \Omega_{a\xi,b\xi}.
  $$
  Since $\xi$ is cyclic, we may approximate any given $\zeta_1$ and $\zeta_2$ by
vectors of the form $a\xi$ and $b\xi$, respectively, and hence one
sees that $\Omega_{\zeta_1,\zeta_2}$ belongs to $K_g$.  This proves
the following:

\state Proposition $K_g$ coincides with
the algebra of all compact operators on $M$.

We now plan to use \lcite{\SufCondForBlend} in order to obtain an affirmative
answer to \lcite{\MyQuestion}.  Under \lcite{\StandHypCommut} we have
that $D$ contains the unit of $A$, so the first part of the hypothesis
of \lcite{\SufCondForBlend} is granted.
  Therefore   we must only verify that
  $e\lambda(A)f$
  consists of compact operators in order to reach our conclusion.

\state Theorem \label HilbScmidt For every $a\in A$, one has that $e\lambda(a)f$ is a
Hilbert-Schmidt operator on $M$.  Moreover, denoting
the Hilbert-Schmidt norm by $\|\cdot\|_2$, one has that
  $$
  \|e\lambda(a)f\|_2\leq \|E(a^*a)\|\half
  \for a\in A.
  $$
  As a consequence $e\lambda(A)f$ consists of compact operators and
hence $\big(K_e,K_f,i_e,i_f,K_g\big)$ is a C*-blend.

\Proof
  Notice that $e\lambda(a)f$ vanishes on the orthogonal complement of $f(M)$.
It is therefore enough to prove that
  $$
  \sum_{i\in I}\|eaf\eta_i\|^2 \leq \|E(a^*a)\|,
  $$
  for any (and hence all) orthonormal basis $\{\eta_i\}_{i\in I}$ of
$f(M)$.
  Since the left hand side is defined to be the supremum of the sums
over finite subsets of $I$, it is enough to prove that
  $$
  \sum_{i=1}^n\|eaf\eta_i\|^2 \leq \|E(a^*a)\|,
  \subeqmark FineteSumForHS
  $$
  for any finite orthonormal set $\{\eta_i\}_{i=1}^n$ contained in the
range of $f$.  By \lcite{\RangeOfF} one has that $f(M)$ is the closure
of $C\xi$ in $M$ so, for each $i$ we may write
  $$
  \ds\eta_i = \lim_{k\to\infty} c_i^k\xi,
  $$
  where $\{c_i^k\}_{k\in\N}$ is a sequence in $C$. Let $\mu^k$ be
the $n\times n$ scalar matrix with $\mu^k_{ij} =
\<c_i^k\xi,c_j^k\xi\>$, as in \lcite{\MainIneq}.
 For all $k$ one then has that
  $$
  \soma \|eac_i^k\xi\|^2 \leq
  \|E(a^*a)\|\ \|\mu^k\|.
  $$
  Since the $\eta_i$ form an orthonormal set, we have that $\mu^k$
converges to the identity matrix, so \lcite{\FineteSumForHS} follows
by taking the limit as $k\to\infty$.
  \endProof

  \endgroup


\Address
  {Departamento de Matem\'atica, Universidade Federal de Santa Catarina,
88040-970 Florian\'opolis SC, Brazil.}
  {exel@mtm.ufsc.br}

\bye